\documentclass[10pt,twoside]{article}
\usepackage{amssymb}
\usepackage[mathscr]{eucal}
\usepackage{eufrak}
\usepackage{amsmath}
\usepackage{mathrsfs}
\font\bg=cmbx10 scaled 1200

\def\intl#1{\int\limits_{#1}}
\def\intll#1#2{\int\limits_{#1}^{#2}}

\def\chiu{\hfill$\displaystyle\vspace{4pt}
\underset\Box\null$}
\def\Pr{{\bf Proof. }}
\def\O{\Omega}

\def\R{\Bbb R}
\def\N{\Bbb N}
\def\o{\"{o}}
\def\à{\`{a}}
\def\è{\`{e}}
\def\ì{\`{i}}
\def\ù{\`{u}}
\def\ò{\`{o}}
\def\é{\'{e}}
\def\vf{\varphi}
\def\dy{\displaystyle}
\def\ve{\varepsilon}
\def\l{\lambda}
\def\pa{\partial}
\def\be{\begin{equation}}
\def\ba{\begin{array}}
\def\ea{\end{array}}
\def\ee{\end{equation}}
\def\vs1{\vspace{1ex}}
\def\cD{{\mathcal D}}
\def\vp{\varphi}
\def\ov{\overline}
\def\po{\partial\Omega}
\font\sc=cmcsc10
\def\intl#1{\int\limits_{#1}}
\def\intll#1#2{\int\limits_{#1}^{#2}}
\title{\bg On a parabolic p-Laplacian system 
with a convective term}
\author{\sc Francesca Crispo and Angelica Pia Di Feola
\thanks{
Dipartimento di Matematica e Fisica, 
Universit\`{a} degli Studi della Campania ``Luigi Vanvitelli'', viale Lincoln 5, 81100 Caserta,
 Italy.
francesca.crispo@unicampania.it; angelicapia.difeola@unicampania.it}}
\date{}
\begin{document}
\maketitle
\noindent{\bf Abstract}  - {\small We prove a result of existence of regular solutions and a maximum principle for solutions to a parabolic $p$-Laplacian system with convective term.  
 }\vskip 0.2cm
 \vskip -0.7true cm\noindent
\newcommand{\red}{\protect\bf}
\renewcommand\refname{\centerline
{\red {\normalsize \bf References}}}
\newtheorem{ass}
{\bf Assumption}[section]
\newtheorem{defi}
{\bf Definition}[section]
\newtheorem{tho}
{\bf Theorem}[section]
\newtheorem{rem}
{\sc Remark}[section]
\newtheorem{lemma}
{\bf Lemma}[section]
\newtheorem{coro}
{\bf Corollary}[section]
\newtheorem{prop}
{\bf Proposition}[section]
\renewcommand{\theequation}{\thesection .
\arabic{equation}}
\setcounter{section}{1}
\section*{\normalsize 1. Introduction}
\renewcommand{\theequation}{1.\arabic{equation}}
\renewcommand{\thetho}{1.\arabic{tho}}
\renewcommand{\thedefi}{1.\arabic{defi}}
\renewcommand{\therem}{1.\arabic{rem}}
\renewcommand{\theprop}{1.\arabic{prop}}
\renewcommand{\thelemma}{1.\arabic{lemma}}
\setcounter{equation}{0} \setcounter{lemma}{0} \setcounter{defi}{0}
\setcounter{prop}{0} 
\setcounter{rem}{0} \setcounter{tho}{0}  This
note deals with the following system \be\label{PF}\begin{array}{ll}\dy
u_t-\nabla\cdot\left(|\nabla u|^{p-2}\nabla u\right)\dy= -\delta u\cdot\nabla u\,,&\hskip-0.2cm\textrm{
in }(0,T)\times\O,
\\\dy \hskip2.4cmu(t,x)\dy=0\,, &\hskip-0.2cm \textrm{ on }
(0,T)\times\po,\;\\\dy\hskip2.4cm u(0,x)=u_\circ(x),&\hskip-0.2cm\mbox{ on
}\{0\}\times\O,\end{array}\ee where $\delta\in \R$ is an a-dimensional constant and 
the $p$-growth exponent belongs
to the interval $(\frac 32,2)$. Here we
assume that $\O$ is a bounded
domain of $\R^3$, 
whose
boundary is $C^2$-smooth, 
$u:(0,T)\times\O\to \R^3,
$
 is a vector field, $u_t:=\frac{\pa u}{\pa t}$, $u \cdot\nabla u:=u_k\partial_k\, u$, and we assume 
the data $u_\circ$ in $ L^\infty(\O)$.
\par The previous is a system strictly connected with the 
model for the dynamic of the so called power-law fluids. These are a special class of non-Newtonian fluids. By non-Newtonian fluids are meant those fluids for which the constitutive relation between the Cauchy stress tensor and the shear rate is non-linear. Denoting by $S$ the shear stress tensor, for a Newtonian fluids $S=2\nu \cD u$ holds, where $\nu>0$ is the viscosity and $\cD u$ the symmetric gradient of the velocity $u$. Fluids not obeying to this linear relation are non-Newtonian. Polymeric liquids, biological fluids, among which blood stands out, suspensions are the most common examples of non-Newtonian fluids. Following \cite{Raj}, the main ways of deviation from a Newtonian behavior, which a fluid may possess all at once or only some of, are: the ability of the fluid to shear thin or shear thicken in shear
flows; the presence of non-zero normal stress differences in shear flows; the ability of the fluid to yield stress;
the ability of the fluid to exhibit stress relaxation;
 the ability of the fluid to creep. 
The fluids we are interested in exhibits mainly the first
 among the listed non-Newtonian behaviors.\par 
Assuming that  the viscosity is a function of the shear rate $|\cD (u)|$ (therefore {\it generalized viscosity}), the constitutive relation becomes 
$S(\cD u) = \nu(|\cD u|)\cD u$. Shear thinning involves a decrease in viscosity as function of the shear rate, while shear thickening involves an increase in viscosity as function of the shear rate. Particularly, power-law fluids show the following rheological law: \be\label{gt}
S(\cD u)=2\mu_0(\mu+|\cD u |)^{p-2}\cD u,\ \ \mu_0>0,\, \mu\geq 0,\, p\in (1,\infty),\ee
where $\mu_0$  is a constant whose dimension is compatible with the expression of the tensor. 
Note that a fluid governed by Navier–Stokes equations is the outcome by setting $\mu_0=\nu$ and $p = 2$, that gives a generalized viscosity equal to a constant. Under the assumption \eqref{gt}, the law of conservation of mass and the law of balance of momentum for a power-law incompressible fluid
read as 
\be\label{PLF}\begin{array}{rll}\vs1\dy u_t-\nabla\cdot \left((\mu+|\cD u |)^{p-2}\cD u\right)+
\nabla\pi_u&= - u\cdot\nabla u+f,&\hskip-0.1cm\textrm{ in
}(0,T)\times\O,
\\\dy\vs1\nabla\cdot u&=0, &\hskip-0.1cm\textrm{ in }(0,T)\times \O,
\end{array}\ee 
where  $u$ is the vector field of velocity, $\pi_u$ is the scalar field of pressure, $f$ is the prescribed body force, and we set for simplicity $2\mu_0=1$. 
\par For references, related both to the physical model and to the mathematical
theory of non-Newtonian fluids, we refer to \cite{15, 22, Raj, show1}.
 \par Over the last fifteen years, we have been interested in different aspects of the mathematical analysis of these fluids, in particular regularity of weak solutions, and existence of ``high regular''  solutions. By high regular solution we mean a solution $u\in L^p(0,T; W_0^{1,p}(\O))\cap C([0,T);L^2(\O))$, with $u_t\in L^2(0,T; L^2(\O))$ such that, for any $\ve>0$, $u\in L^\infty(\ve, T; W^{2,q}(\O))$, $u_t\in L^\infty(\ve, T; L^{q}(\O))$, $q>n$. Therefore $u$ and $\nabla u$ can be pointwisely written. \par The study of system \eqref{PF} 
appears as a further strategic step, which is part of a broader study, already undertaken by one of the authors in a series of papers (see, for instance, \cite{CGM1, CMmodS, CMmodNS}), to better understand some issues connected with power-law fluids.  The leading idea of these works has been to {\it encircle} the problem with problems that are similar in several respects, but contain ``alteration'' of the original one. This is the  idea behind papers \cite{CMmodS, CMmodNS}, related to a modified power-law fluid system in the stationary case:\be\label{stokes}
-\nabla \cdot \left((\mu+|\nabla u|^2)^\frac{p-2}{2} \,\nabla u\right) +\nabla\pi_u= -u\cdot\nabla u+f\,,\quad \nabla \cdot u=0\,.\ee
\par  We have to point out that these  changes make the results in \cite{CMmodS, CMmodNS} not interesting from a fluid dynamic point of view, as the corresponding {\it constitutive law } is not compatible with the principle of material
invariance. However, they are functional to better delimit expectations on possible results on the fluid dynamics problem. \par With the same philosophy we have tried to deal with a modified power-law problem also in the evolutionary case. However, it seems that the extension of the arguments in \cite{CMmodS, CMmodNS} cannot be done using similar techniques. Therefore, in order to approach the evolutionary problem, we  took a step back, asking what it was possible to achieve not only replacing the symmetric gradient with the full gradient, but also by completely eliminating the pressure term, see \cite{CMparabolic}. The problem further deviates from the original fluid dynamic one, and becomes a  parabolic $p$-laplacian system. On the other hand,  this system is widely studied in the mathematical community and has its intrinsic interest. In paper \cite{CMparabolic}, high regularity results for solutions to the $p$-Laplacian system are proved, under suitable
restrictions on the exponent $p\in(1,2)$. Further, a
$L^\infty(\O)$-bound for the weak
solution is there proved. \par In light of the goodness of the results in \cite{CMparabolic}, we thought it appropriate  to see what happens by adding the nonlinear character given by the ``convective'' term $u\cdot\nabla u$. 
Actually, even in this case, we are able to prove the existence of a ``sufficiently regular'' solution (see Definition \ref{defnomu} below), together with a maximum principle. 

\par To introduce our results, we start with the following 
  \begin{defi}\label{defnomu}
{\rm Let  $u_{\circ}\in L^2(\O)$. A field
$u\!:(0,T)\times \O\to\R^N$
 is said a weak solution of
system {\rm \eqref{PF}} if  
$ u\in 
L^{p}(0,T; W_0^{1,p}(\O))\cap L^{p'}(\O_T)$ satisfies 
$$\int_0^T\!\!\! [(u,\psi_\tau)-\left(|\nabla u|^{p-2}\,
\nabla u,\nabla \psi\right)-\delta(u\cdot \nabla u, \psi)]\,d\tau=-(u_\circ, \psi(0)),\,   
\forall \psi\in C_0^\infty([0,T)\times \O),$$ and
$$\lim_{t\to 0^+}(u(t), \vp)=(u_\circ, \vp),\ \forall \vp \in C_0^\infty(\O)\,.$$}
\end{defi}\par
 We are in position to state our main results. 
  The first is an existence theorem of {\it regular} solutions corresponding to an initial datum only in $L^\infty(\O)$,   a maximum principle and a continuity result for solutions. 
\begin{tho}\label{existence}{\sl Let $p\in (\frac 32, 2)$. Corresponding to $u_\circ\in L^\infty(\O)$, there exists a weak solution of \eqref{PF} such that, for any $T>0$,
\be\label{P1}\ba{ll}\vs1\dy 
 u\in C([0,T); L^2(\O))\cap L^\infty(\O_T),\\\dy t^\frac{\alpha}{4-p} \nabla u\in L^\infty(0,T; L^2(\O)), t^\frac{\alpha}{2}u\in L^2(0,T; W^{2,\frac{4}{4-p}}(\O)), t^\frac{1}{2} u_t\in
L^{2}(0,T;L^2(\O)),\ea\ee
for any $\alpha>\frac{4-p}{p}$. 
Moreover, the maximum principle holds 
\be\label{mmtp1}\|u(t)\|_\infty\leq \|u_\circ\|_\infty, \forall t\in (0,T).\ee
Finally, $u\in C([0,T); L^q(\O))$, for any $q\in [1,\infty)$, and 
$$\lim_{t\to 0^+}\|u(t)-u_\circ\|_q=0\,.$$ 
 }
\end{tho}
In particular, estimate \eqref{P1}$_2$ says that if the initial datum is just in $L^\infty(\O)$, then $\nabla u$, $D^2 u$ and $u_t$ have a singularity at $t=0$, which is explicitly computed. If the datum is more regular, we can clearly remove the singularity in the origin, and we state the following  result.
\begin{tho}\label{existencer}{\sl Let $p\in (\frac 32, 2)$. Corresponding to $u_\circ\in L^\infty(\O)\cap W_0^{1,2}(\O)$, there exists a weak solution of \eqref{PF} 
such that, for any $T>0$:
\begin{itemize}
\item[i)]
 $\dy \|u\|_{L^{\infty}(0,T;L^2(\O))}+\|\nabla u\|_{L^{p}(\O_T)}\leq c_1\,,$
 \item[ii)] $\|\nabla u\|_{L^{\infty}(0,T;L^2(\O))}+\|D^2 u\|_{L^{2}(0,T; L^\frac{4}{4-p}(\O))}\!+\|u_t\|_{L^2(\O_T)}\leq c_2 ,$
\end{itemize}
with $c_1:=c_1(\|v_\circ\|_\infty)$ and $c_2:=c_2(\|\nabla v_\circ\|_2,\|v_\circ\|_\infty)$.
Moreover, the maximum principle holds 
\be\label{mmtp1r}\|u(t)\|_\infty\leq \|u_\circ\|_\infty, \forall t\in (0,T).\ee
Finally, $u\in C([0,T); L^q(\O))$, for any $q\in [1,\infty)$, and 
$$\lim_{t\to 0^+}\|u(t)-u_\circ\|_q=0\,.$$ 
 }
\end{tho}
\begin{rem}\label{testsmoothN}
We observe that, since $u(t)\in C([0,T);L^2(\O))$, by using a suitable cut-off function in time, we obtain that, for any $0\le s<t< T$
\be\label{testsmoothst}\ba{ll}\dy\int_s^t [(u,\psi_\tau)-\left(|\nabla u|^{p-2}\,
\nabla u,\nabla \psi\right)\,-\delta(u\cdot \nabla u, \psi)]\,d\tau=
(u(t),\psi(t))
-(u(s), \psi(s)),\\ \hfill  \forall \psi\in
C^\infty_0([0,T)\times\O).\ea\ee
Moreover,
resorting to a density argument, we can take the test functions 
$\psi$ in the space $W^{1,2}(0,T;L^2(\O))\cap L^p(0,T;W_0^{1,p}(\O))\cap L^{p'}(\O_T)$, obtaining an equivalent definition of solution.\end{rem}
Further, corresponding to an initial  datum as in Theorem \ref{existencer}, we prove that regular solutions are unique, provided that $p>\frac 53$.
\begin{tho}\label{uni}{\sl Let be $p>\frac 53$. Corresponding to $u_\circ\in L^\infty(\O)\cap W_0^{1,2}(\O)$, there exists a unique solution of \eqref{PF} with the regularity properties stated in Theorem \ref{existencer}. }
\end{tho}
\par The maximum principle in Theorems \ref{existence} and \ref{existencer} extend the results proved in \cite{CMparabolic} for the $p$-Laplacian parabolic system (without a convective term). A maximum principle for the $p$-Laplacian was previously
proved in \cite{DB,DBH, DBUV} only for
equations, and in
\cite{choe} for systems, but locally.
\par We are not aware of results for a system like \eqref{PF}, nor existence and $L^\infty$-bounds for solutions. In connection with the latter problem, we are just aware of the results in \cite{Avrin} and in \cite{MLectures}, where a quasilinear system of the following kind is considered 
$$u_t-\Delta u\dy= g(t,x,u)\cdot \nabla u,$$
with $g(t,x,u)=u$ in \cite{MLectures}, while $g(t,x,u)$ is an arbitrary locally  Lipschitz function in \cite{Avrin}. In both papers, as in the present one, the global existence of the solution is obtained via the $L^\infty$-bound for the solution itself. Our idea highly relies on that in \cite{MLectures},  
where, for the first time to the best of the authors' knowledge, the duality technique has been employed for proving the maximum principle for a quasilinear  system. On the other hand, the idea used in \cite{Avrin} is to employ a comparison with linear equations of the kind
$$v_t - \Delta v= g(x)\cdot \nabla v,$$
where $g$ is a suitable bounded function,  for which a maximum  principle holds. As for the $p$-laplacian system a similar result is not known, clearly we cannot follow a comparison technique.
\par Further, we prove the extinction of the solution under smallness on the initial data.  
\begin{tho}\label{text}{\sl 
Assume that  the initial datum satisfies
\begin{equation}
{|\delta|}\,||u_\circ||^{\frac{2}{p-1}}_\infty||u_\circ||_2^{2-p} < \gamma, 
\label{hyp}
\end{equation} 
with $\gamma \equiv \gamma(\Omega, p).$ Then, for the solution $u$ in Theorem \ref{existence}, there holds $u(t,x) \equiv 0$ for all $t \geq T^\star$, with
$$T^* \leq\frac{p'||u_\circ||_2^{2-p}}{ \left(-{|\delta|}\, ||u_\circ||_\infty^{\frac{2}{p-1}} ||u_\circ||_2^{2-p} + \gamma\right)   (2-p)}\,.$$ }
\end{tho} 
\par In the sequel we will set $\delta=1$, although all arguments hold for any $\delta$. The interest in an arbitrary a-dimensional parameter $\delta$ is in the possibility of a comparison between the result obtained here for system \eqref{PF} and known results for the $p$-laplacian parabolic system, that corresponds to the  case $\delta=0$. We perform this comparison in Theorem \ref{text}, that for $\delta=0$ reproduces the well known extinction result for solutions of the $p$-laplacian system, since no restrictions on the size of the initial datum remains from \eqref{hyp}, and $T^\star \leq\frac{p'||v_0||^{2-p}_2}{(2-p)\gamma }$ according to \cite{DB}. 
\vskip0.2cm {\it Outline of the proof} -  We introduce a  non-singular approximating problem with a Laplacian term, a mollified right-hand side (see \eqref{PFepv}) and a more regular initial datum, for which, by the Galerkin method,  it is possible to obtain the local existence, and, via a duality technique, the maximum principle. These informations enable us to prove global existence for the approximating solution. Subsequently, taking care of realizing estimates independent of the approximating parameters and of the regularity of the initial datum, we deduce that the same global existence result and $L^\infty$-bound hold for a solution of \eqref{PF}, constructed as limit of the approximations. 
 \par Finally, we like to point out that the restriction from below on the exponent $p$ is mainly due to the application of Lemma \ref{LL1}, which allows us to obtain an estimate of weighted second derivatives.  We believe that by giving up estimates for the second derivatives and settling for a weaker solution, it is possible to lower the exponent $p$.\section*{\normalsize 2. Notations and some auxiliary results}
\renewcommand{\theequation}{2.\arabic{equation}}
\renewcommand{\thetho}{2.\arabic{tho}}
\renewcommand{\thedefi}{2.\arabic{defi}}
\renewcommand{\therem}{2.\arabic{rem}}
\renewcommand{\theprop}{2.\arabic{prop}}
\renewcommand{\thelemma}{2.\arabic{lemma}}
\renewcommand{\thecoro}{2.\arabic{coro}}
\setcounter{equation}{0} \setcounter{coro}{0} \setcounter{lemma}{0}
\setcounter{defi}{0} \setcounter{prop}{0} \setcounter{rem}{0}
\setcounter{tho}{0}
The $L^p$-norm is
denoted by $\|\cdot\|_p$ and, if $m\geq0$, the $W^{m,p}$ and
$W_0^{m,p}$-norms are denoted by $\|\cdot\|_{m,p}$. 
For $q\in[1,\infty)$, and $X$ Banach
space with norm $\|\cdot\|_X$, we denote by $L^q(a,b;X)$ the set of all
function $f:(a,b)\to X$ which are measurable and such that 
\mbox{\footnotesize $\dy\intll ab$}$\|
f(\tau)\|^q_Xd\tau= \|f\|_{L^q(a,b;X)}<\infty$. As well as, if
$q=\infty$ we denote by $L^\infty(a,b;X)$ the set of all function
$f:(a,b)\to X$ which are
 measurable and such that ${\rm ess\ sup}_{t\in(a,b)}\ \|
f(t)\|^q_X=\|f\|_{L^\infty(a,b;X)}<\infty$.  We will often use the identity $L^q((0,T); L^q(\O))=L^q(\O_T)$, for $q\in[1,\infty]$, where  $\O_T:=(0,T)\times \O$.
\par
 For  $p\in (1,2)$, $\mu>0$, and a function $v(t,x)\in L^\infty(0,T; W_0^{1,2}(\O))$, 
we set
\be\label{oper3}
a(\mu,v)(t,x):=\left(\mu+|\nabla
v(t,x)|^2\right)^\frac{p-2}{2}\,,\ee
and,  denoting by 
$\widetilde J_{\eta}$ a time-space Friederich's mollifier, we set 
\be\label{amu}
a_\eta(\mu,v)(t,x):=\left(\mu+|\widetilde J_{\eta}(\nabla
v)(t,x)|^2\right)^\frac{p-2}{2}\,.\ee 
For all $t>0$,  we define, a.e. in $s\in (0,t)$,
\be\label{beeta}\ba{ll}\dy\vs1  (B_\eta(s,x))_{i\alpha j\beta}:=
\delta_{ij}\,\delta_{\alpha\beta}\,a_\eta(\mu, v)(t-s,x).
\ea\ee
Further, we set $\hat v=J_\mu(v)$, with $J_\mu$ space Friedrich's mollifier.\par 
Taking into account the previous notation, we analyze semigroup properties for the following
parabolic system 
\be\label{AD1}\begin{array}{ll}\dy\vs1 \vp_s-\nu\Delta \vp- \nabla\cdot
(B_\eta(s,x)\nabla \vp)= \vp\nabla\cdot \hat v+\hat v\cdot \nabla\vp\,,&\hskip-0.2cm\textrm{ in }(0,t)\times\O,
\\\dy\vs1\hskip3.5cm \vp(s,x)=0\,,&\hskip-0.2cm\textrm{ on }
(0,t)\times\po,\;\\\dy\hskip3.5cm \vp(0,x)=\vp_\circ(x),&\hskip-0.2cm\mbox{ on
}\{0\}\times\O.\end{array}\ee 

%
For system (\ref{AD1}) we employ the existence result stated in Lemma \ref{newag1},  for which we refer to for instance, to \cite{LSU}, Theorem
IV.9.1. 
\begin{lemma}\label{newag1}{\sl Assume that $\nu>0$ and let $\vp_\circ(x)\in C_0^{\infty}(\O)$.  Under the above assumptions on the coefficients, there exists
a unique solution $\vp$ of \eqref{AD1}, such that $\vp\in
 L^2(0,t; W^{2,2}(\O)\cap W_0^{1,2}(\O))$, $\vp_s\in
L^2(0, t; L^{2}(\O))$.
 }\end{lemma}

\vskip0.2cm

\begin{lemma}\label{adgt}{\sl
 Assume that $\vp_\circ(x)\in C_0^{\infty}(\O)$
  and let $\vp$ be the unique solution of \eqref{AD1}. Then, there holds 
  \be\label{max11a}\|\vp(s)\|_1
  \leq \|\vp_\circ\|_1,\ \forall s\in [0,t], \ \textrm{uniformly in } \nu>0 \textrm{ and } \eta>0\,.
\ee
 }\end{lemma}
 \Pr
Let us multiply \eqref{AD1}$_1$ by $\vp (\delta+|\vp|^2)^{-\frac
{1}{2}}$, for some $\delta\in (0,1)$, and integrate by parts on $(0,s)\times \O$, with $s\in [0,t]$. 
Then 
\be\label{max13}
\ba{ll}\vs1 \dy
\|(\delta+|\vp(s)|^2)^\frac 12\|_1+\nu \int_0^s\!\int_\O
(\delta+|\vp|^2)^{-\frac {1}{2}}|\nabla\vp|^2dxd\tau\\ \hskip0.3cm \vs1\dy
-\nu\!\int_0^s\!\int_\O (\delta+|\vp|^2)^{-\frac {3}{2}}(\nabla\vp\cdot
\vp)^2dxd\tau \dy + \int_0^s\!\int_\O B_\eta(\tau, x)\,(\delta+|\vp|^2)^{-\frac
{1}{2}}\,|\nabla\vp|^2dxd\tau
\\\dy\hskip0.3cm \vs1  -\!\! \int_0^s\!\int_\O B_\eta(\tau, x)\,(\delta+|\vp|^2)^{-\frac {3}{2}}\,(\nabla\vp\cdot \vp)^2 dxd\tau\\
 =\vs1 \dy
\|(\delta+|\vp(0)|^2)^\frac 12\|_1+
  \int_0^s\!\int_\O\bigg[ (\delta+| \vp|^2)^{-\frac {1}{2}} |\vp|^2\nabla\cdot \hat v+  \hat v\cdot \nabla (\delta+| \vp|^2)^\frac {1}{2} dxd\tau\bigg].\ea\ee Observing that 
  $$(\nu+B_\eta(\tau, x))\bigg[(\delta+|\vp|^2)^{-\frac
{1}{2}}\,|\nabla\vp|^2-(\delta+|\vp|^2)^{-\frac {3}{2}}\,(\nabla\vp\cdot \vp)^2\bigg] \geq 0,$$
and  integrating by parts the last term on the right-hand side of \eqref{max13}, we find
$$\ba{ll}\vs1 \dy 
\|(\delta+|\vp(s)|^2)^\frac 12\|_1\!\leq &\dy  
 \|(\delta+|\vp_0|^2)^\frac 12\|_1\\ &\dy +
  \int_0^s\!\int_\O\bigg[ (\delta+| \vp|^2)^{-\frac {1}{2}} |\vp|^2\nabla\cdot \hat v- \nabla\cdot \hat v (\delta+| \vp|^2)^\frac {1}{2} \bigg]dxd\tau.\ea$$ 
Since  
$$\lim_{\delta\to 0}\bigg[ (\delta+| \vp|^2)^{-\frac {1}{2}} |\vp|^2\nabla\cdot \hat v- \nabla\cdot \hat v (\delta+| \vp|^2)^\frac {1}{2} \bigg]=0,$$
and $$
\bigg| (\delta+| \vp|^2)^{-\frac {1}{2}} |\vp|^2\nabla\cdot \hat v- \nabla\cdot \hat v (\delta+| \vp|^2)^\frac {1}{2} \bigg|\leq |\nabla\cdot \hat v| ,$$
we can apply the Lebesgue dominated convergence theorem, and we get \eqref{max11a}.\chiu

 \begin{lemma}\label{compattificazione}{\sl
 Let $\nabla\psi\in L^2((0,t)\times \O)$, $\nabla v \in L^r((0,t)
 \times\O)$, for some $r>1$, and let $h^m$ be a
 sequence with $\nabla h^{m}$ bounded in $ L^2((0,t)
 \times \O)$, uniformly in $m\in \N$. Then, there exists a subsequence $h^{m_k}$ such
 that
 $$\lim_{k\to \infty}\int_0^s \int_\O\left((\mu+
 |\widetilde J_{\frac {1}{m_k}}(\nabla v)|^2)^\frac{p-2}{2}-
 (\mu+|\nabla v|^2)^\frac{p-2}{2}\right)
\nabla h^{m_k}\cdot  \nabla
\psi \,dx\,d\tau=0\,.$$ }
 \end{lemma}
 \vskip0.2cm\par
In the sequel, we will use the following results, proved in  \cite{CMparabolic}. 
\begin{lemma}\label{LL1}{\sl
Let $\mu>0$. Assume that $v \in
W^{2,2}(\O)\cap W_0^{1,2}(\O)$. Then, for any $\zeta>0$ such that $C_1(\zeta):=\big(\frac{p}{p(p-1)^2-\zeta}\big)^\frac 12$ is a real number, we have 
$$\big\|(\mu+|\nabla v|^2)^\frac{(p-2)}{4} {D^2 v}\big\|_2\,
\leq C_1(\zeta)\big\|(\mu+|\nabla v|^2)^\frac{(p-2)}{4}{\Delta v}
\big\|_2\!\!+\frac{C_2}{\zeta}\left(\|\nabla v\|_p^p+\mu^\frac p2
|\Omega|\right)^\frac 12.$$ 
}
 \end{lemma}\par We also recall some useful inequalities.
 For the following ``reverse'' version of the H\"older inequality we refer, for instance, to \cite{ADAFOU}, Theorem 2.12.
\begin{lemma}\label{reverse}Let $0<q<1$ and $q'=\frac q{q-1}$. If $f\in L^q(\O)$ and $0<\int_\O|g(x)|^{q'}\,dx<\infty$ then
$$\int_\O\left|f(x)g(x)\right|\,dx\ge\left(\int_\O\left|f(x)\right|^q\,dx\right)^{1/q}\left(\int_\O\left|g(x)\right|^{q'}\,dx\right)^{1/q'}.$$
\end{lemma}
\par For the following generalization of Gronwall's inequality we refer to \cite{MPF}, Ch. 12, and references therein.
\begin{lemma}\label{GGI}{\sl Let $\Phi(t)$ be a nonnegative function satisfying the integral inequality
$$\Phi(t)\leq C+\int_{t_0}^t (A(s)\Phi(s)+ B(s)\Phi(s)^\theta)\, ds, \ C\geq 0, \theta \in [0,1), $$
where $A(t)$ and $B(t)$ are continuous nonnegative functions for $t\geq t_0$. Then, we have
\be\label{GGI1}\ba{ll}\dy
\Phi(t)\leq \left\{C^{1-\theta}\exp\bigg[(1-\theta)\!\intll{t_0}t \! A(s)ds\bigg]+ (1-\theta)\intll{t_0}t \!B(s)\exp\bigg[(1-\theta)\intll{s}t \!A(\tau)d\tau\bigg]ds \right\}^{\!\frac{1}{1-\theta}}.
\ea\ee
}
\end{lemma}
\par 
Finally, we recall the following well known compact embedding result, 
frequently used in the sequel, referring for instance, to \cite{show}, Ch. 3.

\begin{lemma}\label{aubin}{\rm [Aubin-Lions]}- {\sl Let $X$, $X_1$, $X_2$ be
Banach spaces. Assume that   $X_1$ is compactly embedded in $X$ and $X$
 is continuously embedded in $X_2$, and that $X_1$ and $X_2$ are reflexive.
 For $1<q,s<\infty$, set $$W=\{\psi \in L^s(0,T; X_1): \psi_t\in L^q(0,T; X_2)\}\,.$$
Then the inclusion $W \subset L^s(0,T; X)$ is compact.} \end{lemma}

\section*{\normalsize 3. Approximating systems}
\renewcommand{\theequation}{3.\arabic{equation}}
\renewcommand{\thetho}{3.\arabic{tho}}
\renewcommand{\thedefi}{3.\arabic{defi}}
\renewcommand{\therem}{3.\arabic{rem}}
\renewcommand{\theprop}{3.\arabic{prop}}
\renewcommand{\thelemma}{3.\arabic{lemma}}
\renewcommand{\thecoro}{3.\arabic{coro}}
\setcounter{equation}{0} \setcounter{coro}{0} \setcounter{lemma}{0}
\setcounter{defi}{0} \setcounter{prop}{0} \setcounter{rem}{0}
\setcounter{tho}{0}
 Let us study the approximating systems
\be\label{PFl}\begin{array}{ll}\vs1\dy
v_t-\nabla\cdot\left((\mu+|\nabla
v|^2)^\frac{p-2}{2}\nabla v\right)= -J_\mu(v)\cdot \nabla v\,,&\hskip-0.2cm\textrm{ in
}(0,T)\times\O,
\\\dy \vs1\hskip3.5cm v(t,x)=0\,,&\hskip-0.2cm\textrm{ on }
(0,T)\times\po,\;\\\dy\hskip3.5cm v(0,x)=v_\circ(x),&\hskip-0.2cm\mbox{ on
}\{0\}\times\O,\end{array}\ee with $\mu>0$, and

\be\label{PFepv}\begin{array}{ll}\vs1\dy v_t-\nu\Delta
v-\nabla\cdot\left((\mu+|\nabla
v|^2)^\frac{p-2}{2}\nabla v\right)= -J_\mu(v)\cdot \nabla v\,,&\hskip-0.2cm\textrm{ in
}(0,T)\times\O,
\\\dy\vs1\hskip4.6cm v(t,x)=0\,,&\hskip-0.2cm \textrm{ on }
(0,T)\times\po,\;\\\dy\hskip4.6cm v(0,x)=v_\circ(x),&\hskip-0.2cm\mbox{ on
}\{0\}\times\O\,.\end{array}\ee %
\par The aim of this section is to prove the existence of solution of system \eqref{PFl} in the following sense.
\begin{defi}\label{wmu} {\rm Let $\mu> 0$. Let $v_{\circ}\in
L^\infty(\O)\cap W_0^{1,2}(\O)$. A field $v\!:(0,T)\times \O\to\R^N$
 is said a solution of
system {\rm \eqref{PFl}} if, for any $T>0$,\begin{itemize}
\item[i)] 
$
 v\in C([0,T); L^2(\O))\cap L^{p}(0,T; W_0^{1,p}(\O))\cap L^\infty(\O_T)$, \par 
 $ \nabla v\in L^\infty(0,T; L^2(\O)), v\in L^2(0,T; W^{2,\frac{4}{4-p}}(\O)), v_t\in
L^{2}(\O_T),$ \par
\item[ii)] $v$ satisfies the integral identity
$$\ba{ll}\dy \vs1\int_0^t \left[(v,\psi_\tau) -\left(a(\mu,
v)\,\nabla v,\nabla \psi\right)-(J_\mu(v)\cdot \nabla v, \psi)\right]d\tau=(v(t),
\psi(t))-(v_\circ, \psi(0)),\\ \hfill   \forall t\in [0,T), \ \forall \psi\in C_0^\infty([0,T)\times \O),\ea$$ 
\item[iii)]$\dy \lim_{t\to
0^+}\|v(t)-v_\circ\|_2=0\,.$\end{itemize}}
\end{defi}
\begin{defi}\label{weakmunu}
{\rm Let $\mu> 0$, $\nu>0$. Let $v_{\circ}\in W_0^{1,2}(\O)$. A field
$v\!:(0,T)\times\O\to\R^N$
 is said a  solution of
system {\rm \eqref{PFepv}} if, for any $T>0,$\begin{itemize}
\item[i)]$ v\in   C([0,T); L^2(\O))\cap L^{\infty}(0,T;
W_0^{1,2}(\O))\cap L^2(0,T; W^{2,2}(\O)), v_t\in   L^2(\O_T)$,\item[ii)]$v$ satisfies the integral identity  $$\ba{ll}\dy\vs1 \int_0^t
\left[(v,\psi_\tau)-\nu(\nabla v,\nabla\psi)-\left(a(\mu, v)\,\nabla
v,\nabla \psi\right)-(J_\mu(v)\cdot \nabla v, \psi)\right]d\tau\\ \hskip4cm \dy=(v(t),\psi(t))-(v_\circ,
\psi(0)),\ \forall t\in [0,T),\
   \forall \psi\in C_0^\infty([0,T)\times \O),\ea$$\item[iii)]$\dy\lim_{t\to
0^+}\|v(t)-v_\circ\|_2=0\,.$\end{itemize}}
\end{defi}


In order to study the existence and regularity of a solution of
 \eqref{PFl}, firstly we study the same issues for
the parabolic approximating system \eqref{PFepv}. The existence
and regularities for the solution of this latter system are obtained 
 in Proposition \ref{existenceL} and \ref{existencemunu}, for a suitably regular datum, and in Corollary \ref{existencepeso} without exploiting the regularity of the initial datum other than the $L^\infty$-property. The idea is to start by proving the local existence result of Proposition \ref{existenceL}, and then obtain, in Proposition \ref{existencemunu}, the global existence throughout the $L^\infty$-estimate of Proposition  \ref{ulinfty}. 

\begin{prop}\label{existenceL}
{\sl Let be $\nu>0$ and $\mu>0$. Assume that $v_\circ$
belongs to $W_0^{1,2}(\O)$. Then there exists $T>0$ and a solution $v$
of system \eqref{PFepv} in [0,T) satisfying {\sl i)--iii)} in  Definition\,\ref{weakmunu}. In particular $T\in (0,\frac{1}{c(\mu,\nu)\|v_\circ\|_2^2})$, where $c(\mu,\nu)$ is a positive constant that blows up as $\nu$ or $\mu$ tend to zero.
}
\end{prop}
  {\bf Proof} -  In the sequel, we prove the existence of a solution by the Galerkin
   method.
   Let $\{a_j\}$ be the eigenfunctions of $-\Delta$, and denote
by $\lambda_j$ the corresponding eigenvalues:
$$\ba{ll}\vs1 \dy-\Delta a_j\!\!\!\!&\dy=\lambda_j a_j,\ \mbox{ in }
\O,\\\dy \hskip0.6cm a_j\!\!\!\!&=\dy0,\hskip0.5cm  \mbox{ on }
\po\,. \ea$$  Recall that $a_j\in
W_0^{1,2}(\O)\cap W^{2,2}(\O)$, it is a basis in $W_0^{1,2}(\O)$, and it is orthonormal in $L^2(\O)$, with 
 $(\nabla a_j, \nabla a_k)=( a_j, a_k)=0$ for $j\not= k$,
 $\|\nabla a_j\|_2^2=\lambda_j$. \par
 We consider the Galerkin
approximations related to system \eqref{PFepv} of the form
\be\label{G1}v^k(t,x)=\sum_{j=1}^k c_{jk}(t)a_j(x)\ \hskip 0.2cm
k\in \N\,,\ee where the coefficients $c_{jk}$ satisfy the following system
of ordinary differential equations
\be\label{lpsc}\ba{ll}\vspace{0.5ex} \dy \dot
c_{jk}(t)=-\nu\sum_{i=1}^k b_{ji}c_{ik}(t)-
\sum_{i=1}^k d_{ji}c_{ik}(t)-\sum_{i=1}^{k}\sum_{l=1}^kB_{jil}c_{ik}(t)c_{lk}(t), \hskip0.2cm  j=1,\ldots,k,\\
\dy c_{jk}(0):=(v_{\circ},a_j),\ea\ee with $b_{ji}:=(\nabla a_i,\nabla
a_j)$, $d_{ji}:=((\mu+|\nabla (c_{rk}(t)
a_r)|^2)^\frac{(p-2)}{2}\,\nabla a_i,\nabla a_j)$, $B_{jil} :=(J_{\mu}(a_i)\cdot\nabla a_l, a_j)$, $i, l=1,\cdots, k$. With this choice
of $c_{jk}(t)$ we impose that $v^{k}(t,x)$ are solutions of the
following system of $k$-differential equations \be\label{ap1}
(v^k_t,a_j)-\nu(\Delta v^k, a_j)+(\nabla\cdot 
\left(a(\mu,v^k)\nabla v^k
\right)
, a_j)+(J_\mu(v^k)\cdot\nabla v^k, a_j)=0\,,\ 
j=1,\cdots,k,\ee with initial conditions $c_{jk}(0)=(v_\circ, a_j)$,
$j=1,\cdots, k$. 
We explicitly remark that $(\nabla v_\circ, \nabla a_j)=\lambda_j(v_\circ, a_j)$, hence $\nabla v^k(0,x)=\sum_{j=1}^kc_{jk}(0)\nabla a_j= \sum_{j=1}^k (\nabla v_\circ,\nabla a_{j})\frac{\nabla a_j}{\lambda_j}$ that gives
$$\|\nabla v^k(0)\|_2^2=\sum_{j=1}^k(\nabla v_\circ,\frac{\nabla a_{j}}{\|\nabla a_j\|_2})^2\leq \|\nabla v_\circ\|_2^2.$$
\par As the right-hand side of \eqref{lpsc} is a locally 
Lipschitz function, due to the assumption $\mu>0$, 
the existence of a solution to \eqref{ap1} in a maximal time interval
$[0,T_k]$ follows by standard results on ordinary
differential equations.  The following a priori estimate will ensure that there exists a suitable constant $c$ such that $T_k\geq  \frac{1}{c\|v_\circ\|_2^2}$, for all $k\in \N$. 
\vskip0.2cm {\underline {\it Energy estimate} - Let us multiply
\eqref{ap1} by $c_{jk}$
and sum over $j$. We get
 \be\label{ap2}
\frac12\frac{d}{dt}\|v^k\|_2^2+\nu\|\nabla v^k\|_2^2+\|(\mu+|\nabla
v^k|^2)^\frac{(p-2)}{4}\nabla v^k\|_2^2=\|(J_\mu(v^k)\cdot\nabla v^k)\cdot v^k\|_1\,.\ee
By H\o lder's inequality and using that $\|J_\mu(v^k)\|_\infty\leq c(\mu)\|v^k\|_2$,  with $c(\mu)\to \infty$ as $\mu\to 0$, 
the {\it
energy identity} \eqref{ap2}
implies 
$$\frac12\frac{d}{dt}\|v^k\|_2^2+\nu\|\nabla v^k\|_2^2+\|(\mu+|\nabla
v^k|^2)^\frac{(p-2)}{4}\nabla v^k\|_2^2\leq c(\mu)\|\nabla v^k\|_2\|v^k\|_2^2\,.$$
By applying Young's inequality on the right-hand side, one gets
\be\label{ap2abb}
\frac12\frac{d}{dt}\|v^k\|_2^2+\frac{\nu}{2}\|\nabla v^k\|_2^2+\|(\mu+|\nabla
v^k|^2)^\frac{(p-2)}{4}\nabla v^k\|_2^2\leq \frac{2c(\mu)}{\nu}\|v^k\|_2^4=:c(\mu, \nu)\|v^k\|_2^4
.\ee
Further, by using, on the left-hand side, the inequalities
\be\label{gp}\ba{ll}\vs1\dy\int_{\O}|\nabla v^k|^p\, dx\!\!\dy
=\int_{|\nabla v^k|^2> \mu}|\nabla v^k|^p\, dx+\int_{|\nabla v^k|^2
\leq \mu}|\nabla v^k |^p\, dx\\ \dy\leq 2^\frac{2-p}{2}\!\int_{\O}
a(\mu,v^k)\,|\nabla v^k|^2 dx+\int_{|\nabla v^k|^2\leq \mu}\!\!\mu^\frac
p2\, dx\\ \dy \leq  2^\frac{2-p}{2} \int_{\O} a(\mu,v^k)\, |\nabla v^k|^2dx+ \mu^\frac p2\,
|\O|,\ea\ee
one gets
\be\label{ap2ab}
\frac12\frac{d}{dt}\|v^k\|_2^2+\frac{\nu}{2}\|\nabla v^k\|_2^2+2^{\frac{p-2}{2}}\|\nabla
v^k\|_p^p\leq c(\mu, \nu)\|v^k\|_2^4+2^{\frac{p-2}{2}}\mu^\frac p2\,
|\O|
\,.\ee
We end up with the following differential inequality
\be\label{ap2abc}
\frac12\frac{d}{dt}(\|v^k\|_2^2+c(\mu,\nu,p,\O))\leq c(\mu,\nu)(\|v^k\|_2^2+c(\mu, \nu, p,\O))^2\,.\ee
 The integration of the above differential inequality, yields 
\be\label{ap2ac}
\|v^k(t)\|_2^2\leq \frac{\|v_\circ\|_2^2}{1-c(\mu,\nu)\|v_\circ\|_2^2  t},\ \forall t\in\big[0, \frac{1}{c(\mu,\nu)\|v_\circ\|_2^2}\big),\ee
%
that, for $T<\frac{1}{c(\mu,\nu)\|v_\circ\|_2^2}$, gives
\be\label{gb}\|v^k(t)\|_2^2=|c_k(t)|^2\leq
\frac{\|v_\circ\|_2^2}{1-c(\mu,\nu)\|v_\circ\|_2^2  t},\ \forall \ t\in [0, T]\subseteq [0,T_k)\,, \ee
for all $k\in \N$. 
\par
Once (\ref{gb}) is at disposal, an integration in time of estimate \eqref{ap2abb} and estimate \eqref{ap2ab} ensures that
\be\label{nunu}\|\nabla
v^k\|^2_{L^2(\O_T)}\leq c(\mu,\nu)\,(1+\|v_\circ\|_2^4)\,.\ee
\vskip0.2cm {\underline {\it Further regularity} -  Let us multiply \eqref{ap1} by $\l_j c_{jk}$ , and   
sum
over $j$.  Then, calculating the divergence of $ 
a(\mu,v^k)\nabla v^k$, and using H\o lder's inequality, we have
\be\label{inc1}\ba{ll}\vs1\dy
 \frac12\frac{d}{dt}\|\nabla
v^k\|_2^2+\nu\|\Delta v^k\|_2^2+ \,\|(\mu+|\nabla
v^k|^2)^\frac{(p-2)}{4}\Delta v^k\|_2^2\\ \dy \vs1\leq
(2-p) \,\|(\mu+|\nabla
v^k|^2)^\frac{(p-2)}{2}|D^2v^k||\Delta v^k|\|_1 +\|J_\mu(v^k)\|_\infty\|\nabla v^k\|_{2}\|\Delta v^k\|_{2} \,.\ea\ee
For the first term on the right-hand side, we employ H\o lder's inequality and then Lemma\,\ref{LL1}, and we find, for any $\zeta$ satisfying the assumptions of Lemma \ref{LL1},  
\be\label{inc2}
\ba{ll}\vs1\dy\,\|(\mu+|\nabla
v^k|^2)^\frac{(p-2)}{2}|D^2v^k||\Delta v^k|\|_1 \leq 
\,\|(\mu+|\nabla
v^k|^2)^\frac{(p-2)}{4}D^2v^k\|_2\,\|(\mu+|\nabla
v^k|^2)^\frac{(p-2)}{4}\Delta v^k\|_2\\\vs1
\dy \leq \left[
C_1(\zeta)\,\|(\mu+|\nabla
v^k|^2)^\frac{(p-2)}{4}\Delta v^k\|_2+ \frac{C_2}{\zeta}\,\left(
\|\nabla v^k\|_p^p+ \mu^\frac p2|\O|\right)^\frac 12\right]
\|(\mu+|\nabla
v^k|^2)^\frac{(p-2)}{4}\Delta v^k\|_2\\
\dy =
C_1(\zeta)\,\|(\mu+|\nabla
v^k|^2)^\frac{(p-2)}{4}\Delta v^k\|_2^2+ \frac{C_2}{\zeta}\,\left(
\|\nabla v^k\|_p^p+ \mu^\frac p2|\O|\right)^\frac 12\|(\mu+|\nabla
v^k|^2)^\frac{(p-2)}{4}\Delta v^k\|_2.
\ea
\ee
We insert this estimate in \eqref{inc1}. Then we apply Cauchy-Schwartz inequality with any $\delta>0$, and we get \be\label{ap5n}\ba{ll}\vs1\dy
 \frac12\frac{d}{dt}\|\nabla
v^k\|_2^2+\nu\|\Delta v^k\|_2^2+ \,\|(\mu+|\nabla
v^k|^2)^\frac{(p-2)}{4}\Delta v^k\|_2^2\\ \dy \vs1\leq
\left((2-p)\,C_1(\zeta)+\frac{\delta}{2}\right) \,\|(\mu+|\nabla
v^k|^2)^\frac{(p-2)}{4}\Delta v^k\|_2^2+ \frac{C}{2\delta}\,\left(
\|\nabla v^k\|_p^p+ \mu^\frac p2|\O|\right)\\ \dy\hfill +\|J_\mu(v^k)\|_\infty\|\nabla v^k\|_{2}\|\Delta v^k\|_{2}.\ea\ee
Let us fix $\zeta=\frac{p(2p-3)}{2}$. Since $p>\frac 32$,
we can  choose $\delta=1-(2-p)C_1$ in such a way that $\overline C(p):=1-(2-p)C_1-\frac \delta2>0 $. By applying on the last term Cauchy-Schwartz inequality, we finally get 
\be\label{appri6n}\ba{ll}\vs1\dy \frac 12\frac{d}{dt}\|\nabla
v^k\|_2^2+\frac \nu 2 \|\Delta v^k\|_2^2+\overline C(p)
\,\|(\mu+|\nabla
v^k|^2)^\frac{(p-2)}{4}\Delta
v^k\|_2^2\\ \dy \hskip2cm\leq C \|\nabla
v^k\|_p^p+\, C\,\mu^\frac p2
|\O|+c(\nu)\|J_\mu(v^k)\|_\infty^2\|\nabla v^k\|_{2}^2\,. \ea\ee 
Analogously, multiplying \eqref{ap1} by $d c_{jk}/dt$, 
summing 
over $j$, and then applying on the right-hand side H\o lder's and Cauchy inequality, we find
 \be\label{ap4an}\ba{ll}\vs1\dy  \|v^k_t\|_2^2+\frac
{\nu}{2}\frac{d}{dt}\|\nabla v^k\|_2^2+\frac
1p\frac{d}{dt}\|(\mu+|\nabla
v^k|^2)^\frac{1}{2}\|_p^p=-(J_\mu(v^k)\cdot \nabla v^k, v^k_t)\\ \hfill\dy
\leq \|J_\mu(v^k)\|_\infty\|\nabla v^k\|_2\|v^k_t\|_2\leq \frac 12 \|v^k_t\|_2^2+\frac 12\|J_\mu(v^k)\|_\infty^2\|\nabla v^k\|_2^2.\ea\ee %
Summing \eqref{appri6n} and \eqref{ap4an}, then integrating between $0$ and $t$, we find 
\be\label{apppri67n}\ba{ll}\vs1\dy (1+\nu)\|\nabla v^k(t)\|_2^2+\nu 
\int_0^t\|\Delta v^k\|_2^2d\tau+
\int_0^t\|v^k_\tau\|_2^2d\tau+\frac
2p\|(\mu+|\nabla
v^k(t)|^2)^\frac{1}{2}\|_p^p
\\ \vs1\dy\hfill
\leq (1+\nu)\|\nabla v_\circ\|_2^2+\frac
2p\|(\mu+|\nabla
v_\circ|^2)^\frac{1}{2}\|_p^p\\ \hfill \dy+\,C\int_0^t (\|\nabla
v^k\|_p^p+\,\mu^\frac p2
|\O|)d\tau   +c(\mu,\nu)\int_0^t \|v^k\|_2^2\|\nabla v^k	\|_2^{2}d\tau. 
\ea\ee  As the time integrals on the right-hand side can be bounded in terms of $\|v_\circ\|_2$ using the energy estimate, the previous inequality ensures that 
$$\{\nabla v^k\}\ \mbox{ is bounded in } L^\infty(0,T; L^2(\O)), \mbox{uniformly
with respect to  }k,$$ 
$$\{D^2 v^k\}\mbox{ is bounded in } L^2(0,T; L^2(\O)), \mbox{uniformly
with respect to  }k,$$ 
%
%
and  $$\{v^k_t\}
\ \mbox{ is bounded in } L^2(0,T; L^2(\O)), \mbox{uniformly
with respect to  }k.$$\par

\vskip0.2cm {\underline {\it Limiting process} -    
Using the boundedness of the sequence $\{v^k\}$ in the spaces listed above, we can extract a subsequence, still
denoted by $\{v^k\}$, such that, in the limit as $k$ tends to
$\infty$,
$$ v^k\, \rightharpoonup \, v
\textrm{ in } L^{\infty}(0,T; L^2(\O)) \ \textrm{weakly}-*\,;$$
 $$ v^k\, \rightharpoonup \, v \textrm{ in }
L^{2}(0,T; W^{2,2}(\O))\ \textrm{weakly}\,,$$
$$ v^k\, \rightarrow \, v \textrm{ in }
L^{2}(0,T; W^{1,r}_0(\O))\ \textrm{strongly}\,, \mbox{for any} \ r<6,$$
$$v^k_t\rightharpoonup v_t\
\textrm{in } L^{2}(0,T; L^2(\O))\ \textrm{weakly}.$$
As far as the convective term, writing for all $\psi\in C_0^\infty([0,T)\times \O)\, $
$$\ba{ll}\vs1 \dy|\int_0^t(J_\mu(v^k)\cdot \nabla v^k, \psi)d\tau-\int_0^t (J_\mu(v)\cdot \nabla v, \psi)d\tau|\\\vs1 \dy= |\int_0^t(J_\mu(v^k-v)\cdot \nabla v^k, \psi)d\tau+\int_0^t(J_\mu(v)\cdot (\nabla v^k-\nabla v), \psi)d\tau|
\ea$$ 
using that $\|J_\mu(v^k-v)\|_\infty\leq c(\mu)\|v^k-v\|_2$, and then the strong convergence of $v^k$ to $v$ in $L^2(0,T; L^2(\O))$ and the weak convergence of $\nabla v^k$ to $\nabla v$ in $L^2(\O_T)$, one shows that the above terms tend to zero. 
 Following usual arguments (see proof of Proposition 3.4 for details), one
shows that the limit $v$ is a weak solution of
system \eqref{PFepv}. The regularities stated for $v$ follows from
the analogous regularities of $\{v^k\}$ and the lower semi-continuity of
the norm for the weak convergence.\chiu

\begin{prop}\label{ulinfty} {\sl 
Let $v$ be the solution of \eqref{PFepv} corresponding to
$v_\circ\in L^\infty(\O)\cap W_0^{1,2}(\O)$. Then \be\label{infv}\|v(t)\|_\infty\leq
\|v_\circ\|_\infty, \ \mbox{ for all }\ t\in [0,T).\ee}
\end{prop}
\Pr  Let us consider, for a fixed $\eta>0$,  the solution $\vp^\eta(s,x)$
of system \eqref{AD1} corresponding to a data $\vp_\circ\in
C_0^\infty(\O)$, where $B_\eta$ is given by \eqref{beeta} and $\hat v=J_\mu(v)$,
 and set
$\hat\vp^\eta(\tau)=\vp^\eta(t-\tau)$. 
Then, using $\hat\vp^\eta(\tau)$ as test
function in the weak formulation of \eqref{PFepv}, we
have
 \be\label{a8}\ba{ll}\dy\vs1 (v(t),\vp_\circ)-(v_\circ,\vp^\eta(t))
-\int_0^t (v(\tau),\hat\vp_\tau^\eta(\tau))d\tau
\\\dy\hskip2cm
=-\nu\int_0^t(\nabla v(\tau),\nabla
\hat\vp^\eta(\tau))d\tau-\int_0^t(a(\mu,v(\tau)) \nabla
v(\tau),\nabla \hat\vp^\eta(\tau))d\tau\\
\dy \hskip2.4cm -\int_0^t (J_\mu(v)\cdot \nabla v, \hat\vp^\eta(\tau))d\tau.\ea\ee
We write the second 
term on the right-hand side of \eqref{a8}, as
\be\label{mts2}\ba{ll}\dy \vs1\int_0^t(a(\mu,v)
\nabla v(\tau),\nabla \hat\vp^\eta(\tau))d\tau=\!\!&\dy\! \int_0^t(
\nabla v(\tau),a_\eta(\mu,v(\tau))\nabla
\hat\vp^\eta(\tau))d\tau\\&\dy +\int_0^t( \nabla v(\tau),\nabla
\hat\vp^\eta(\tau))[a(\mu,v) -a_\eta(\mu, v)]\,d\tau\,.\ea
\ee 
Further, we integrate by parts the third term on the right-hand side of \eqref{a8}, and we get
\be\label{mts1}\ba{ll}\dy \int_0^t\!\! (J_\mu(v)\cdot\! \nabla v, \hat\vp^\eta(\tau))d\tau=-\!
\!\int_0^t\! (\nabla\!\cdot J_\mu(v) v, \hat\vp^\eta(\tau))d\tau-\int_0^t\! (v,J_\mu(v)\!\cdot\nabla\hat\vp^\eta(\tau))d\tau
\ea\ee
Since $\vp^\eta$ is the solution of \eqref{AD1}, taking into account \eqref{mts1} and \eqref{mts2}, identity
\eqref{a8} becomes \be\label{a9}\ba{ll}\dy\vs1
(v(t),\vp_\circ)=(v_\circ,\vp^\eta(t))\\ \dy \hskip0.5cm -\int_0^t( \nabla
v(\tau),\nabla \hat\vp^\eta(\tau))[a(\mu, v(\tau))-a_\eta(\mu,
v(\tau))]\,d\tau 
=:(v_\circ,\vp^\eta(t))+I_\eta.\ea\ee
By applying Lemma\,\ref{compattificazione} with $\nabla h^\eta=\nabla \hat\vp^\eta\in
L^2(0,t;L^2(\O))$ and with $\nabla\psi= \nabla v\in L^{2}(0,t; L^2(\O))$,
due to  Lemma\,\ref{newag1} and to Proposition\,\ref{existence}, 
the integral
$I_\eta$ goes to zero, as $\eta$ goes to zero, along a subsequence.
Finally, using \eqref{max11a} and then passing to the limit as $\eta$
tends to zero  in \eqref{a9}, along a subsequence, we get
$$|(v(t),\vp_\circ)| \leq \|v_\circ\|_{\infty}\,\|\vp_\circ\|_{1}, \
\forall \vp_\circ\in C_0^\infty(\O),$$
and, using density arguments, for any $\varphi_\circ\in L^1(\O)$.  This last  implies $$
\|v(t)\|_{\infty}=\!\!\sup_{\vf_\circ\in L^1(\O) \atop
|\vf_\circ|_1=1}|(v(t),\vf_\circ)| \leq
\|v_\circ\|_{\infty}\,.$$ 
\chiu

\begin{prop}\label{existencemunu}
{\sl Let $\nu>0$ and $\mu>0$. Assume that $v_\circ$
belongs to $L^\infty(\O)\cap W_0^{1,2}(\O)$. Then there exists a solution $v$
of system \eqref{PFepv} satisfying, for all $T>0$, the following bounds, \underline{uniformly in $\nu>0$, $\mu>0$}:\begin{itemize}
\item[i)]
$\dy \|v\|_{L^\infty(0,T;L^2(\O))}+\|\nabla v\|_{L^{p}(\O_T)}\leq K_1(\|v_\circ\|_\infty,T)\,;$
 \item[ii)] $\|\nabla v\|_{L^{\infty}(0,T;L^2(\O))}+\|D^2 v\|_{L^{2}(0,T; L^\frac{4}{4-p}(\O))}+\|v_t\|_{L^2(\O_T)}\leq\,K_3 (\|\nabla v_\circ\|_2,\|v_\circ\|_\infty,T)\,.$
\end{itemize}}
\end{prop}
\Pr From Proposition\,\ref{existenceL},  there exists a function $v$,
solution of \eqref{PFepv} in the time interval $[0, T)\subset[0, \frac{1}{c(\mu,\nu)\|v_\circ\|_2^2})$ corresponding to the initial data
$v_\circ$. As, by assumption,  $v_\circ\in L^\infty(\O)$, from Proposition \ref{ulinfty} we also know that $v\in L^{\infty}(\O_T)$ and, for any $t\in[0,T)$ satisfies the stimate
\be\label{v2i}
\|v(t)\|_2\leq |\O|^\frac12 \|v_\circ\|_\infty.\ee
The validity of the above inequality ensures that the solution $v$ actually exists for any $T>0$.
\par
Let us use $v$ as test function in the Definition \ref{weakmunu}, by density,  and integrate by parts. We get
 \be\label{app2}\ba{ll}\dy\vs1
\|v(t)\|_2^2+2\nu \int_0^t\|\nabla v(\tau)\|_2^2\,d\tau+2\int_0^t\|(\mu+|\nabla
v(\tau)|^2)^\frac{(p-2)}{4}\nabla v(\tau)\|_2^2\,d\tau\\\hfill \dy=\|v_\circ\|_2^2-\int_0^t(J_\mu(v)\cdot\nabla v, v)\,d\tau\,.\ea\ee
By H\o lder's inequality, the property $\|J_\mu(v)\|_\infty\leq \|v\|_\infty\leq \|v_\circ\|_\infty$, and then Young's and convexity inequality, 
we find
$$|(J_\mu(v)\cdot\nabla v, v)|\leq \|J_\mu(v)\|_\infty \|\nabla v\|_p\|v\|_{p'}\leq  \frac 1p\|\nabla v\|_p^p+\frac{1}{p'}\|v_\circ\|_\infty^\frac{2}{p-1} \|v\|_2^2.$$
Inserting the above estimate in 
identity \eqref{app2}, and recalling estimate \eqref{gp},
yields  
 \be\label{app2a}
\ba{ll}\dy\vs1
\|v(t)\|_2^2+2\int_0^t\nu\|\nabla v(\tau)\|_2^2d\tau+c\int_0^t\|\nabla v(\tau)\|_p^pd\tau\\\dy\leq \|v_\circ\|_2^2+ \mu^\frac p2\, |\O|t
+\frac{1}{p'}\|v_\circ\|_\infty^\frac{2}{p-1}\int_0^t  \|v(\tau)\|_2^2d\tau. \ea\ee
From \eqref{v2i} and \eqref{app2a} we obtain
\be\label{gpunu}\|v\|_{L^\infty(0,T; L^2(\O))}+\|\nabla v\|_{L^p(\O_T)}\leq K_1(\|v_\circ\|_\infty, T),\ee
with a constant $K_1$ independent of $\nu$ and $\mu$, that 
also implies
  $$\|(\mu+|\nabla
v|^2)^\frac{(p-2)}{2}\nabla
v\|^{p'}_{L^{p'}(\O_T)}\leq \|\nabla
v\|^p_{L^{p}(\O_T)}\leq \,
K_1\,.$$ 
\vskip0.2cm {\underline {\it Further regularity} -  Let us  multiply the equation \eqref{PFepv} by $-\Delta v$\footnote{The difference with the approach used in the proof of Proposition \ref{existenceL} is that now we search for  estimates that are uniform in $\nu$ and $\mu$.}.
By employing H\o lder's inequality, Lemma\,\ref{LL1} with $\zeta=\frac{p(2p-3)}{2}$,  and Cauchy-Schwartz inequality with any $\delta>0$, we get \footnote{For details on this computation, see \eqref{inc1}--\eqref{ap5n}}\be\label{ap5}\ba{ll}\vs1\dy
 \frac12\frac{d}{dt}\|\nabla
v\|_2^2+\nu\|\Delta v\|_2^2+ \,\|(\mu+|\nabla
v|^2)^\frac{(p-2)}{4}\Delta v\|_2^2\\ \dy \vs1\leq
\left((2-p)\,C_1(\zeta)+\frac{\delta}{2}\right) \,\|(\mu+|\nabla
v|^2)^\frac{(p-2)}{4}\Delta v\|_2^2+ \frac{C}{2\delta}\,\left(
\|\nabla v\|_p^p+ \mu^\frac p2|\O|\right)\\ \dy\hfill +\|J_\mu(v)\|_\infty\|\nabla v\|_{\frac4p}\|\Delta v\|_{\frac{4}{4-p}}
 \,,\ea\ee
 for any $\delta>0$. Hence, since $p>\frac 32$,
 choosing $\delta=1-(2-p)C_1$ in such a way that $\overline C(p):=1-(2-p)C_1-\frac \delta2>0 $,  we get 
\be\label{appri6}\ba{ll}\vs1\dy \frac 12\frac{d}{dt}\|\nabla
v\|_2^2+\nu \|\Delta v\|_2^2+\overline C(p)
\,\|(\mu+|\nabla
v|^2)^\frac{(p-2)}{4}\Delta
v\|_2^2\\ \dy \hskip2cm\leq c \|\nabla
v\|_p^p+\, c\,\mu^\frac p2
|\O|+\|v_\circ\|_\infty\|\nabla v\|_{\frac4p}\|\Delta v\|_{\frac{4}{4-p}}\,. \ea\ee On the left-hand side we apply the reverse H\o lder inequality of Lemma \ref{reverse}, with exponents $q=\frac{2}{4-p} \in (0,1)$, $q'=\frac{2}{p-2}<0$, and we get
$$\|(\mu+|\nabla
v|^2)^\frac{(p-2)}{4}\Delta
v\|_2^2\geq (\mu|\O|+\|\nabla v\|_2^2)^\frac{p-2}{2}\|\Delta v\|^2_{\frac{4}{4-p}}=:B(\mu,\nabla v)^\frac{p-2}{2}\|\Delta v\|^2_{\frac{4}{4-p}}. $$ By using this inequality in \eqref{appri6}, observing that
\be\label{star}B(\mu,\nabla v)^\frac{2-p}{2}\frac 12\frac{d}{dt}\|\nabla
v\|_2^2=\frac 1{4-p}\frac{d}{dt}B(\mu,\nabla v)^\frac{4-p}{2}\, \ee
we find
\be\label{flu2}\ba{ll}\vs1\dy \frac 1{4-p}\frac{d}{dt}B(\mu,\nabla v)^\frac{4-p}{2}+\overline C(p)
\,\|\Delta
v\|_{\frac{4}{4-p}}^2\\ \dy \hskip0.2cm\leq c B(\mu,\nabla v)^\frac{2-p}{2}(\|\nabla
v\|_p^p+\, \,\mu^\frac p2
|\O|)+B(\mu,\nabla v)^\frac{2-p}{2}\|v_\circ\|_\infty\|\nabla v\|_{\frac4p}\|\Delta v\|_{\frac{4}{4-p}}\,. \ea\ee 
On the last term on the right-hand side we apply Gagliardo-Nirenberg
inequality, and we get
$$ \|\nabla v\|_\frac4p\leq  c\|D^2 v\|^a_\frac{4}{4-p}\|v\|_q^{1-a}\leq c\|D^2 v\|^a_\frac{4}{4-p}\|v\|_\infty^{1-a},$$ with $\dy a=\frac{12+4q-3pq}{12-4q+3pq}\,\in \left(\frac 12,1 \right)\ \mbox{for } q<\frac{4}{3p-4}.$ 
Hence, \eqref{flu2} becomes
$$\ba{ll}\vs1\dy \frac 1{4-p}\frac{d}{dt}B(\mu,\nabla v)^\frac{4-p}{2}+\overline C(p)
\,\|\Delta
v\|_{\frac{4}{4-p}}^2\\ \dy \hskip0.2cm\leq c B(\mu,\nabla v)^\frac{2-p}{2}(\|\nabla
v\|_p^p+\, \,\mu^\frac p2
|\O|)+cB(\mu,\nabla v)^\frac{2-p}{2}\|v_\circ\|_\infty^{2-a}\|\Delta v\|^{1+a}_{\frac{4}{4-p}}\,, \ea$$
where we have also used $\|D^2 v\|_\frac{4}{4-p}\leq c\|\Delta  v\|_\frac{4}{4-p}$. 
Young's inequality with exponents $\frac{2}{1+a}$ and $\frac{2}{1-a}$ now gives
\be\label{apppri6}\ba{ll}\vs1\dy \frac 1{4-p}\frac{d}{dt}B(\mu,\nabla v)^\frac{4-p}{2}+(\overline C(p)-\ve)
\,\|\Delta
v\|_{\frac{4}{4-p}}^2\\ \dy \hskip0.2cm\leq c B(\mu,\nabla v)^\frac{2-p}{2}(\|\nabla
v\|_p^p+\, \,\mu^\frac p2
|\O|)+c(\ve) B(\mu,\nabla v)^{\frac{2-p}{2}\frac{2}{1-a}}\|v_\circ\|_\infty^{\frac{2(2-a)}{1-a}}\,. \ea\ee
Dividing both side by $B(\mu,\nabla v)^\frac{2-p}{2}$, and using \eqref{star}, we arrive at the following differential inequality
$$\ba{ll}\vs1\dy \frac 1{2}\frac{d}{dt}\|\nabla v\|_2^2\leq c (\|\nabla
v\|_p^p+\,\mu^\frac p2
|\O|)+c B(\mu,\nabla v)^{\frac{2-p}{2}(\frac{2}{1-a}-1)}\|v_\circ\|_\infty^{\frac{2(2-a)}{1-a}}\,. \ea$$
By integrating between $0$ and $t$, recalling the definition of $B$, and setting $$\theta:=\frac{2-p}{2}(\frac{2}{1-a}-1),$$we find 
\be\label{apppri67}\ba{ll}\vs1\dy \|\nabla v(t)\|_2^2\leq \|\nabla v_\circ\|_2^2+c\int_0^t (\|\nabla
v\|_p^p+\,\mu^\frac p2
|\O|)d\tau \\ \dy\hskip2.3cm + \int_0^t (\mu|\O|)^\theta \|v_\circ\|_\infty^{\frac{2(2-a)}{1-a}}d\tau +\int_0^t \|v_\circ\|_\infty^{\frac{2(2-a)}{1-a}}\|\nabla v	\|_2^{2\theta}d\tau \,. \ea\ee It is easy to check that, by choosing $q>\frac{12-6p}{3p-4}$ too, there holds  $\theta<1.$ 
Hence, using the nonlinear Gronwall's inequality of Lemma \ref{GGI} and estimate {\sl i)},  yields
\be\label{appppri6}\ba{ll}\vs1\dy \|\nabla v(t)\|_2^2\leq  K_2(\|\nabla v_\circ\|_2,\|v_\circ\|_\infty),\ea\ee
with $K_2$ independent on $\mu\leq 1$, but depending on $\|v_\circ\|_{1,2}$. Using estimate \eqref{appppri6} in \eqref{apppri6} ensures that $$D^2 v\in L^2(0,T; L^\frac{4}{4-p}(\O)), \mbox{uniformly
with respect to  }\nu, \mu,$$
 and 
\be\label{sd1}\ba{ll}\vs1\dy \|D^2 v(t)\|_{L^2(0,T; L^\frac{4}{4-p}(\O))}\leq  K_3(\|\nabla v_\circ\|_2,\|v_\circ\|_\infty).\ea\ee
Now, let us multiply  \eqref{PFepv} by $v_t$, and then apply on the right-hand side H\o lder's and Cauchy's inequalities. We find
 \be\label{ap4a}\ba{ll}\vs1\dy  \|v_t\|_2^2+\frac
{\nu}{2}\frac{d}{dt}\|\nabla v\|_2^2+\frac
1p\frac{d}{dt}\|(\mu+|\nabla
v|^2)^\frac{1}{2}\|_p^p=-(J_\mu(v)\cdot \nabla v, v_t)\\ \hfill\dy
\leq \|J_\mu(v)\|_\infty\|\nabla v\|_2\|v_t\|_2\leq \frac 12 \|v_t\|_2^2+\frac 12\|v_\circ\|_\infty^2\|\nabla v\|_2^2.\ea\ee %
Integrating between $0$ and $t$, and using estimate \eqref{appppri6} we find
$$\ba{ll}\dy \int_0^t\|v_\tau\|_2^2d\tau+\nu \|\nabla v(t)\|_2^2+\frac
2p\|(\mu+|\nabla
v(t)|^2)^\frac{1}{2}\|_p^p\\ \hfill\dy \leq \nu \|\nabla v_\circ\|_2^2+\frac
2p\|(\mu+|\nabla
v_\circ|^2)^\frac{1}{2}\|_p^p+\|v_\circ\|_\infty^2\int_0^t\|\nabla v(\tau)\|_2^2d\tau,\ea$$
that ensures $v_t\in L^{2}(\O_T)$.\par
Finally, the continuity property of the solution $v$ in $L^2(\O)$ follows from $v\in$ \\ $ L^{p}(0,T; W_0^{1,p}(\O))$ coupled with  $v_t\in L^{2}(\O_T)$.\chiu
\par
If we do not exploit the $W_0^{1,2}$-
regularity of the initial data, but just its $L^\infty$-property, we
have the following corollary.
\begin{coro}\label{existencepeso}
{\sl Under the assumptions of Proposition \ref{existencemunu}, there exists a constant $K$, independent of $\|\nabla v_\circ\|_2$ but depending on $\|v_\circ\|_\infty$,  such that 
\begin{itemize}
 \item[ii')] $\|t^{\frac{\alpha}{4-p}}\nabla v\|_{L^{\infty}(0,T;L^2(\O))}+\|t^{\frac{\alpha}{2}}D^2 v\|_{L^{2}(0,T; L^\frac{4}{4-p}(\O))}+\|t^\frac 12\,v_t\|_{L^2(\O_T)}\leq K,$
\end{itemize}
\par where $\alpha>\frac{4-p}{p}$. }
\end{coro}
\Pr
We multiply equation \eqref{PFepv} by $-\Delta v$. As in the proof of Proposition \ref{existencemunu}, after some manipulations, we find estimate \eqref{flu2},
that for the readers' convenience we reproduce here
\be\label{fnu}\ba{ll}\vs1\dy \frac 1{4-p}\frac{d}{dt}B(\mu,\nabla v)^\frac{4-p}{2}+\overline C(p)
\,\|\Delta
v\|_{\frac{4}{4-p}}^2\\ \dy \hskip0.2cm\leq C B(\mu,\nabla v)^\frac{2-p}{2}(\|\nabla
v\|_p^p+\, \,\mu^\frac p2
|\O|)+B(\mu,\nabla v)^\frac{2-p}{2}\|v_\circ\|_\infty\|\nabla v\|_{\frac4p}\|\Delta v\|_{\frac{4}{4-p}}\,, \ea\ee 
with $$B(\mu,\nabla v):=(\mu|\O|+\|\nabla v\|_2^2).$$
For $t>0$, let us multiply inequality \eqref{fnu} by $t^\alpha$ with $\alpha>0$ to be chosen:
\be\label{III}\begin{array}{ccc}\vs1
\dfrac{d}{dt}\left( t^\alpha B(\mu,\nabla v)^\frac{4-p}{2} \right) +t^\alpha \overline C(p) ||\Delta v||_\frac{4}{4-p}^2 \leq 
 c t^\alpha \left( ||\nabla v||_p^p + \mu^\frac{p}{2}
|\Omega|\right) B(\mu,\nabla v)^\frac{2-p}{2}  \\\vs1
\hspace*{1.5cm} \,\,\,+ \,\alpha t^{\alpha -1} B(\mu,\nabla v)^\frac{4-p}{2} 
+\,  t^\alpha ||v_\circ||_\infty ||\nabla v||_\frac{4}{p} ||\Delta v ||_{\frac{4}{4-p}} B(\mu,\nabla v)^\frac{2-p}{2}\\\hspace*{1.5cm} =: I_1 + I_2 + I_3.
\end{array}\ee
We leave $I_1$ as it is. To estimate $I_2$, we apply Gagliardo-Nirenberg inequality and we have  
$$ ||\nabla v||_2^{4-p} \leq c ||D^2v||_{\frac{4}{4-p}}^{\beta(4-p)} ||v||_{\overline{q}} ^{(1-\beta)(4-p)},$$
with 
$$\frac{1}{2} = \frac{1}{3} + \beta \left( \frac{4-p}{4} - \frac{2}{3}\right) + (1-\beta) \frac{1}{\overline{q}}.$$
So,
$$\beta(\ov q):=\beta = \dfrac{2(6-\overline{q})}{3p\overline{q} +12 -4 \overline{q}} \in \left( \frac{1}{2} , 1 \right)\ \mbox{iff} \ 
\overline{q} < \frac{4}{p}.$$
Moreover, since 
\begin{equation}
 \frac{2}{\beta(4-p)}>1,
\label{condq}
\end{equation} we can apply Young's inequality and we have, for any $\varepsilon>0$,\\
$$\begin{array}{ccc}
I_2=\alpha t^{\alpha-1} B(\mu, \nabla v)^\frac{4-p}{2} \leq c t^{\alpha-1}(\mu|\O|)^\frac{4-p}{2}+\varepsilon \left( \alpha t^{\frac{\alpha}{2}\beta(4-p)}||\Delta v||^{\beta(4-p)}_{\frac{4}{4-p}}\right)^\frac{2}{\beta (4-p)} \; \\ \hspace*{4cm}+ \;  c(\varepsilon) \left(\alpha t^{\alpha-1-\frac{\alpha}{2}\beta(4-p)} ||v||_{\overline{q}}^{(1-\beta)(4-p)} \right)^{\frac{2}{2-\beta (4-p)}},
\end{array}$$
where we have also used $\|D^2 v\|_\frac{4}{4-p}\leq c\|\Delta  v\|_\frac{4}{4-p}$.
Similarly, for $I_3$, by using Gagliardo-Nirenberg inequality we have
$$||\nabla v||_\frac{4}{p} \leq c ||D^2v||_{\frac{4}{4-p}}^a ||v||_{\overline{\overline{q}}}^{1-a},$$
with
$$\frac{p}{4} = \frac{1}{3} + a \left(\frac{4-p}{4} - \frac{2}{3} \right) + \left( 1-a\right) \frac{1}{{\overline{\overline{q}}}},$$
where, since $p>\frac 32$, we have  \be\label{aovq}a(\ov{\ov q}):=a=\frac{12+4 {\overline{\overline{q}}} - 3p {\overline{\overline{q}}}}{12 + 3 p {\overline{\overline{q}}} -4{\overline{\overline{q}}}}\in (\frac 12, 1)\mbox{ iff }\ {\overline{\overline{q}}} < \frac{4}{3p-4}.\ee
Therefore, applying Young's inequality with exponent $\frac{2}{1+a},$ we find 
$$\ba{ll}\vs1\dy 
I_3 \leq  c t^{\alpha} ||v_\circ||_\infty  B(\mu, \nabla v)^\frac{2-p}{2} ||v||_{{\overline{\overline{q}}}}^{1-a}||\Delta v||_{\frac{4}{4-p}}^{1+a}
\\\vs1 \dy \leq \varepsilon\left( t^{\alpha\frac{1+a}{2}} ||D^2v||^{a+1}_{\frac{4}{4-p}} \right)^\frac{2}{1+a} + c(\varepsilon) \left( t^{\alpha - \alpha (\frac{1+a}{2})} ||v_\circ||_\infty B(\mu, \nabla v)^\frac{2-p}{2} ||v||_{{\overline{\overline{q}}}}^{1-a}\right)^\frac{2}{1-a}. 
\end{array}$$
 Substituting the previous estimates in \eqref{III}, we finally get 
 $$\begin{array}{ll}\vs1
 \dfrac{d}{dt}( t^\alpha B(\mu, \nabla v)^\frac{4-p}{2})+  t^\alpha (\overline C(p)-\varepsilon) ||\Delta v||^2_{\frac{4}{4-p}}\\ \dy \hfill \leq c t^\alpha \left( ||\nabla v||_p^p + \mu^\frac{p}{2} |\Omega| \right) B(\mu, \nabla v)^\frac{2-p}{2}  +\, ct^\alpha ||v_\circ||_{\infty}^{\frac{2}{1-a}} B(\mu, \nabla v)^\frac{2-p}{1-a}||v||_{\overline{\overline{q}}}^2 
  \\ \hspace{2.75cm} + \,  c \alpha\left(t^{\alpha -1-\frac{\alpha}{2}\beta (4-p)} ||v||_{\ov q}^{(1-\beta)(4-p)}\right)^{\frac{2}{2-\beta(4-p)}}+c t^{\alpha-1}(\mu|\O|)^\frac{4-p}{2}.
 \end{array}$$
By defining  \be\label{phi}
\Phi(t) := t^\alpha B(\mu, \nabla v)^\frac{4-p}{2},\ee and estimating the norms $\|v\|_{\ov q}$ 
and $\|v\|_{\ov{\ov q}}$ by $\|v_\circ\|_\infty$, we rewrite the previous inequality as follows  
$$\begin{array}{ccc}\vs1
\dfrac{d \Phi}{dt}   + t^\alpha (\overline C(p)-\varepsilon) ||\Delta v||_\frac{4}{4-p}^2\leq c\,\Phi^{\frac{2-p}{4-p}} \left( ||\nabla v||_p^p + \mu^\frac{p}{2} |\Omega|\right) t^{\alpha \frac{2}{4-p}} \,+
\\ \vs1\dy\hspace{4.2cm} +\, c\,\Phi^{\frac{2(2-p)}{(4-p)(1-a)}} 
t^{\alpha \left( 1- \frac{2(2-p)}{(4-p)(1-a)}\right)} ||v_\circ||_{\infty}^{\frac{2}{1-a} +2} +\\
 \hspace{2.65cm} + \; c \alpha t^{\alpha - \frac{2}{2-\beta(4-p)}} ||v_\circ||_{\infty}^{\frac{2(1-\beta)(4-p)}{2-\beta(4-p)}}+c t^{\alpha-1}(\mu|\O|)^\frac{4-p}{2}.
\end{array}$$
By recalling \eqref{aovq}, we require $$\frac{6(2-p)}{3p-4}<\overline{\overline{q}}<\frac{4}{3p-4}.$$ This ensures that 
 $\frac{2(2-p)}{(4-p)(1-a)}<1$. 
Since $a \in \left( \frac{1}{2}, 1\right)$, setting
\be\label{teta}
 \theta:=\max \lbrace \frac{2-p}{4-p}, \; \frac{2}{1-a}\frac{2-p}{4-p} \rbrace =\frac{2}{1-a}\frac{2-p}{4-p}\,,\ee we have
$$\begin{array}{l}\vs1
\dfrac{d\Phi}{dt} + t^\alpha (\overline C(p)-\varepsilon) ||\Delta v||_\frac{4}{4-p}^2\leq c\left( 1 + \Phi^\theta \right) t^{\alpha \frac{2}{4-p}} \left(||\nabla v||_p^p+ \mu^\frac{p}{2} |\Omega|\right) \\+ c\Phi^\theta  t^{\alpha \left( 1 - \frac{2(2-p)}{(4-p)(1-a)} \right)}||v_\circ||_\infty^\frac{4-2a}{1-a}
+\; ct^{\alpha - \frac{2}{2-\beta (4-p)}} ||v_\circ||_\infty^{\frac{2(1-\beta)(4-p)}{2-\beta(4-p)}}+c t^{\alpha-1}(\mu|\O|)^\frac{4-p}{2}.
\end{array}$$
Integrating on $(0,t)$ we obtain
$$\begin{array}{ll}\vs1
\displaystyle \Phi(t) +\int_0^t \tau^\alpha (\overline C(p)-\varepsilon) ||\Delta v||_\frac{4}{4-p}^2 d\tau \\\vs1 \dy \leq  c\int_0^t\!\!\left( 1 + \Phi^\theta(\tau) \right) \tau^{\alpha \frac{2}{4-p}} \left(||\nabla v||_p^p + \mu^\frac{p}{2}|\Omega| \right) d\tau 
 \displaystyle +c\int_0^t\!\Phi^\theta(\tau) \tau^{\alpha\left( 1-\frac{2(2-p)}{(4-p)(1-a)}\right)}||v_\circ||_\infty^{\frac{4-2a}{1-a}} d\tau\\\dy +c\int_0^t \tau^{\alpha-\frac{2}{2-\beta(4-p)}} ||v_\circ||_\infty^{\frac{2(4-p)(1-\beta)}{2-\beta(4-p)}}d\tau+c\int_0^t  \tau^{\alpha-1}(\mu|\O|)^\frac{4-p}{2}
 \, d\tau.\end{array}$$
We choose $\alpha$ in such a way that \be\label{cal}
\alpha - \frac{2}{2-\beta(4-p)}> -1.\ee
We search for the smallest exponent $\alpha$ satisfying \eqref{cal}. Recalling that, from condition \eqref{condq}, $\frac12<\beta<\frac{2}{4-p}$, and setting $f(\beta):=-1+\frac{2}{2-\beta(4-p)}$, this function is increasing in $\beta\in [\frac 12, \frac{2}{4-p})$ and $f(\beta)>0$, which means $\frac{2}{2-\beta(4-p)} >1$, since it is the conjugate exponent of $\frac{2}{\beta(4-p)}$. 
Therefore $f(\beta)$ reaches its minimum value in $\beta=\frac 12$ and we find that \eqref{cal} is satisfied as soon as $\alpha>f(\frac 12)=\frac{4-p}{p}$. \par
We can now use Lemma \ref{GGI}, with the following choice:
$$C:= cT^{\alpha \frac{2}{4-p}} \int_0^T \left(||\nabla v||_p^p + \mu^\frac{p}{2}|\Omega| \right) d\tau  +c T^{\alpha-\frac{2}{2-\beta(4-p)}+1} ||v_\circ||_\infty^{\frac{2(4-p)(1-\beta)}{2-\beta(4-p)}}+cT^\alpha \mu^\frac{1}{2}|\Omega|,  $$
$$A(s) =0,$$ 
$$\ba{ll}\vs1B(s):\!\!\!\!&\dy = s^{\alpha \frac{2}{4-p}} \left(||\nabla v(s)||_p^p + \mu^\frac{p}{2}|\Omega| \right) + s^{\alpha\left( 1-\frac{2(2-p)}{(4-p)(1-a)}\right)}||v_\circ||_\infty^{\frac{4-2a}{1-a}}\\&\dy = s^{\alpha \frac{2}{4-p}} \left(||\nabla v(s)||_p^p + \mu^\frac{p}{2}|\Omega| \right) + s^{\alpha( 1-\theta)}||v_\circ||_\infty^{\frac{4-2a}{1-a}}\;\ea $$
By the bounds given in estimate {\sl i)} of Proposition \ref{existencemunu} and recalling the definition \eqref{phi} of $\Phi$, we  get
$$t^{\alpha} \|\nabla v\|_2^{4-p}\leq C+(T^{\alpha \frac{2}{4-p}}\,K_1^p(\|v_\circ\|_\infty))^\frac{1}{1-\theta}+ K_4(\|v_\circ\|_\infty)(t^{(\alpha \frac{2}{4-p}+1)\frac{1}{(1-\theta)}}+ t^{(\alpha( 1-\theta)+1)\frac{1}{(1-\theta)}}),$$
with clear meaning of constant $K_4$. 
Dividing by $t^\alpha$ and rasing to the power $\frac{1}{4-p}$,  observing that 
${(\alpha \frac{2}{4-p}+1)\frac{1}{(1-\theta)}}>\alpha$ and 
$ (\alpha( 1-\theta)+1)\frac{1}{(1-\theta)}>{\alpha}$, 
we can simply write $$\|\nabla v(t)\|_2\leq \frac{c}{t^\frac{\alpha}{4-p}}+K_4(\|v_\circ\|_\infty)t^\gamma,$$
for a suitable $\gamma>0$. Therefore
we have $t^{\frac{\alpha}{4-p}} \nabla v$ and 
  $t^\frac\alpha2 D^2v $ are bounded in $
  L^\infty(0, T; L^2(\Omega))$ and in $L^2(0,T; L^\frac{4}{4-p}(\Omega))$, respectively, by a constant depending on the initial datum only throughout 
 $ \|v_\circ\|_\infty$. 
\\
 About the estimate on the time derivative, let us start from  inequality \eqref{ap4a}:
 $$||v_t||_2^2 + \dfrac{\nu}{2} \dfrac{d}{dt} \|\nabla v \|_2^2 + \dfrac{1}{p} \dfrac{d}{dt} \left| \left| \left( \mu + |\nabla v|^2 \right)^\frac{1}{2} \right| \right|_p^p \leq ||v_\circ||_\infty ||\nabla v||_2 ||v_t||_2.$$
 Multiplying by t and applying  and Cauchy-Schwartz inequality we find
$$\begin{array}{l}\vs1
t||v_t||_2^2 + \dfrac{\nu}{2} \dfrac{d}{dt} \left(t \|  \nabla v \|_2^2 \right)+ \dfrac{1}{p} \dfrac{d}{dt} \left(t\left| \left| \left( \mu + |\nabla v|^2 \right)^\frac{1}{2} \right| \right|_p^p \right)\\ 
\hspace{3.7cm}  \leq t ||v_\circ||_\infty ||\nabla v||_2 ||v_t||_2 \; + \,\dfrac{\nu}{2} ||\nabla v||_2^2\, + \,\dfrac{1}{p} \left| \left| \left( \mu + |\nabla v|^2\right)^\frac{1}{2}\right| \right|_p^p  \\
\hspace{3.7cm} \leq \dfrac{t}{2} ||v_\circ||^2_\infty ||\nabla v||^2_2 + \dfrac{t}{2} ||v_t||^2_2 +\dfrac{\nu}{2} ||\nabla v||_2^2\, + \,\dfrac{1}{p} \left| \left| \left( \mu + |\nabla v|^2\right)^\frac{1}{2}\right| \right|_p^p. 
\end{array}$$
Integration between $0$ and $t$ yields
$$\begin{array}{l}
\vs1
\displaystyle\int_0^t\tau||v_\tau||_2^2d\tau + {\nu} t\, \|  \nabla v(t) \|_2^2 + \dfrac{2}{p}  t\,\left| \left| \left( \mu + |\nabla v(t)|^2 \right)^\frac{1}{2} \right| \right|_p^p  \\ 
\hspace{2.5cm}\displaystyle \leq\int_0^t{\tau} ||v_\circ||^2_\infty ||\nabla v||^2_2 d\tau \displaystyle + {\nu}\int_0^t  ||\nabla v||_2^2 d\tau \, + \,\dfrac{2}{p} \int_0^t\left| \left| \left( \mu + |\nabla v|^2\right)^\frac{1}{2}\right| \right|_p^p d\tau.
\end{array}$$
By the uniform estimates of the right-hand side, $t^\frac{1}{2}v_t \in L^2(\Omega_T)$, with a norm bounded by a constant depending on the initial datum only throughout 
 $ \|v_\circ\|_\infty$. 
 \chiu

\begin{prop}\label{existencemu}[Convergence of the $\nu$-solutions] 
{\sl Let $\mu>0$. Assume that $v_\circ$
belongs to $L^\infty(\O)\cap W_0^{1,2}(\O)$. Then there exists a solution $v$
of system \eqref{PFl} in the sense of Definition\,\ref{wmu}.
In particular, for any $T>0$:\begin{itemize}
\item[i)]
$\dy \|v\|_{L^{\infty}(0,T;L^2(\O))}+\|\nabla v\|_{L^{p}(\O_T)}\leq K_1(\|v_\circ\|_\infty,T)\,;$
 \item[ii)] $\|\nabla v\|_{L^{\infty}(0,T;L^2(\O))}+\|D^2 v\|_{L^{2}(0,T; L^\frac{4}{4-p}(\O))}+\|v_t\|_{L^2(\O_T)}\leq\,K_3 (\|\nabla v_\circ\|_2,\|v_\circ\|_\infty,T).$
\end{itemize}
Moreover, $v$ satisfies the maximum principle
\be\label{mpvmu}\|v(t)\|_\infty\leq \|v_\circ\|_\infty, \ \forall t\in (0,T).\ee}
\end{prop}

\Pr For any fixed $\nu>0$, denote by $v^\nu$ the solution of \eqref{PFepv} corresponding to the initial data
$v_\circ\in W_0^{1,2}(\O)\cap L^\infty(\O)$, satisfying    the
bounds {\sl i), ii)} in Proposition\,\ref{existencemunu},  uniformly in $\nu>0$, $\mu>0$.
Hence we can extract a
subsequence, still denoted by
$\{v^\nu\}$, satisfying the following convergence properties:
\be\label{wnoee} v^\nu\, \rightharpoonup \, v
\textrm{ in } L^{\infty}(0,T; L^2(\O))\ \textrm{weakly}^*\,,\ee
\be\label{wnoe} v^\nu\, \rightharpoonup \, v
\textrm{ in } L^{p}(0,T; W^{1,p}_0(\O))\ \textrm{weakly}\,,\ee
\be\label{weg} \nabla v^\nu\, \rightharpoonup \, \nabla v
\textrm{ in } L^{\infty}(0,T; L^2(\O))\ \textrm{weakly}^*\,,\ee
\be\label{weds} v^\nu\, \rightharpoonup \, v \textrm{ in }
L^{2}(0,T; W^{2, \frac{4}{4-p}}(\O))\ \textrm{weakly}\,,\ee
\be\label{wedt}v^\nu_t\rightharpoonup v_t\
\textrm{in } L^{2}(0,T; L^{2}(\O))\ \textrm{weakly}.\ee
By applying Lemma \ref{aubin}, from the last two convergences we also find
\be\label{dn1} v^\nu\, \rightarrow \, v
\textrm{ in } L^{2}(0,T; W^{1,r}_0(\O))\ \textrm{strongly}\,,\mbox{ for any } r<(\frac{4}{4-p})^*\,,\ee
where, for any $q>1$, by $q^*$ we denote the Sobolev conjugate exponent of $q$, i.e. $q^*=\frac{3q}{3-q}$.  
The above strong convergence implies both the almost everywhere convergence of a (not relabeled) subsequence of $\{\nabla v^\nu\}$ to $\nabla v$, and the following strong convegence \be\label{dn2} v^\nu\, \rightarrow \, v
\textrm{ in } L^{2}(\O_T).\ee
Further, by standard arguments, one can verify that, for any $\vp\in L^2(\O)$, $$(v^\nu, \vp) \mbox{ is equibounded and equicontinuous in } t\in [0,T),$$
whence,  from the Ascoli-Arzel\`a theorem and the weak  convergence of the sequence $\{v^\nu\}$ to $v$ in $L^2(\O_T)$, we get, up to a subsequence, 
\be\label{wut}\lim_{\nu\to 0^+}(v^\nu(t)-v(t), \vp)=0,\ \mbox{ uniformly in } t\in [0,T)\mbox{ and for all } \vp\in L^2(\O),\ee
and $(v(t), \vp)\in C([0,T)), \forall \vp\in L^2(\O).$
\par Further, by elementary calculations and H\o lder's inequality we get
\be\label{ado}\ba{ll}\vs1\dy \|v^{\nu_1}(t)-v^{\nu_2}(t)\|_2^2=\int_0^t\frac{d}{d\tau}\|v^{\nu_1}(\tau)-v^{\nu_2}(\tau)\|_2^2d\tau\\\vs1\hfill\dy =2
\int_0^t\,(v^{\nu_1}(\tau)-v^{\nu_2}(\tau))\cdot (v^{\nu_1}_\tau(\tau)-v_\tau^{\nu_2}(\tau))d\tau\\ \hfill\dy\leq
2\|v^{\nu_1}-v^{\nu_2}\|_{L^2(\O_t)}\|v_\tau^{\nu_1}-v_\tau^{\nu_2}\|_{L^2(\O_t)}.
\ea \ee
By the strong convergence of the sequence $\{v^{\nu}\}$ in $L^2(0,T; L^2(\O))$ and the boundedness of $ v^{\nu}_t$ uniformly in $\nu$ in $L^2(0,T; L^2(\O))$, we find that for any fixed $t\geq 0$ the sequence $\{v^\nu(t)\}$ is a Cauchy sequence in $L^2(\O)$. Therefore, for any fixed $t\geq 0$ the sequence $\{v^\nu(t)\}$ strongly converges in $L^2(\O)$, whence a.e. in $\O$. We can write, for any $t\geq 0$ and a.e. in $\O$, 
$$|v(t,x)|\leq |v^\nu(t,x)-v(t,x)|+\|v^\nu(t)\|_\infty\leq |v^\nu(t,x)-v(t,x)|+\|v_\circ\|_\infty,$$
where, in the last step, we have used Proposition \ref{ulinfty} on $v^\nu$. 
Passing to the limit on $\nu$, then taking the essential supremum, we also get the validity of \eqref{mpvmu} for $v$. \par
Thus, everything is arranged to pass to the limit in the weak formulation of the
approximations,   with test function  in $C_0^\infty([0,T)\times \O)$:
$$\ba{ll}\dy\vs1 \int_0^t
\left[(v^\nu,\psi_\tau)-\nu(\nabla v^\nu,\nabla\psi)-\left(a(\mu, v^\nu)\,\nabla
v^\nu,\nabla \psi\right)-(J_\mu(v^\nu)\cdot \nabla v^\nu, \psi)\right]d\tau\\ \hskip0.5cm \dy -(v^\nu(t),\psi(t))+(v_\circ,
\psi(0))=:\sum_{i=1}^5 I_i^\nu+(v_\circ,
\psi(0)).
\ea$$
Indeed 
$$I_1^\nu \to \int_0^t
(v,\psi_\tau)d\tau $$
from the strong convergence of $v^\nu$ in $L^2(\O_T)$; as, from the weak convergence of $\nabla v^\nu$ in $L^2(0,T; L^r(\O))$ to $\nabla v$, we have
$$\int_0^t
(\nabla v^\nu,\nabla\psi)d\tau \to\int_0^t
(\nabla v,\nabla\psi)d\tau ,$$
then $$
I_2^\nu\to 0.$$
 We also have  that
$$a(\mu, v^\nu)\,\nabla v^\nu \mbox{ boundend in } L^{p'}(\O_T) $$ and, due to the already exploited almost everywhere convergence of $\{\nabla v^\nu\}$ to $\nabla v$, we have  $$a(\mu, v^\nu)\,\nabla v^\nu \rightarrow a(\mu, v)\,\nabla v, \mbox { a.e. in } \ \O_T.$$ 
Then, Lemma I.1.3 in \cite{lions}
 allows us to
infer the validity of the weak convergence
$$a(\mu, v^\nu)\,\nabla v^\nu \rightharpoonup a(\mu, v)\,\nabla v,  \mbox{ weakly in } L^{p'}(\O_T), $$ 
 that ensures 
 $$I_3^\nu\to \int_0^t (a(\mu, v)\,\nabla v, \nabla \psi) d\tau.$$
 As far as the convective term is concerned, by using H\o lder's inequality we have 
 $$\ba{ll}\vs1 \dy|I_4^\nu\!-\!\!\int_0^t\!(J_\mu(v)\cdot \nabla v, \psi)d\tau|= |\int_0^t\!(J_\mu(v^\nu-v)\cdot \nabla v^\nu, \psi)d\tau\!+\!\int_0^t\!(J_\mu(v)\cdot (\nabla v^\nu-\nabla v), \psi)d\tau|\\ \dy \leq\|v^\nu-v\|_{L^2(\O_T)}\|\nabla v^\nu\|_{L^2(\O_T)}\|\psi\|_{L^\infty(\O_T)}+\|v\|_{L^2(\O_T)}\|\nabla v^\nu-\nabla v\|_{L^2(\O_T)}\|\psi\|_{L^\infty(\O_T)}.  \ea$$ 
 By the strong convergence of $v^\nu$ to $v$ and of $\nabla v^\nu$ to $\nabla v$, both in $L^2(\O_T)$, the above difference tends to zero. 
Finally, by employing estimate \eqref{wut}, we have  $I_5^\nu\to (v(t), \psi(t)). $
\par Properties {\sl i), ii)}
 follow
from the analogous of
$v^\nu$ in Proposition \ref{existencemunu} and the lower
semi-continuity of the norm for the
weak convergence.\chiu
 \begin{coro}\label{existencemup}{\sl 
Under the assumptions of Proposition \ref{existencemu} there exists a constant $K$, independent of $\|\nabla v_\circ\|_2$ but depending on $\|v_\circ\|_\infty$,  such that \begin{itemize}
 \item[ii')] $\|t^{\frac{\alpha}{4-p}}\nabla v\|_{L^{\infty}(0,T;L^2(\O))}+\|t^{\frac{\alpha}{2}}D^2 v\|_{L^{2}(0,T; L^\frac{4}{4-p}(\O))}+\|t^\frac 12\,v_t\|_{L^2(\O_T)}\leq K,$
\end{itemize}
\par where $\alpha>\frac{4-p}{p}$. }
\end{coro}

\vskip0.5cm%
\section*{\normalsize 4. Proof of the main theorems}
\renewcommand{\theequation}{4.\arabic{equation}}
\renewcommand{\thetho}{4.\arabic{tho}}
\renewcommand{\thedefi}{4.\arabic{defi}}
\renewcommand{\therem}{4.\arabic{rem}}
\renewcommand{\theprop}{4.\arabic{prop}}
\renewcommand{\thelemma}{4.\arabic{lemma}}
\renewcommand{\thecoro}{4.\arabic{coro}}
\setcounter{equation}{0} \setcounter{coro}{0} \setcounter{lemma}{0}
\setcounter{defi}{0} \setcounter{prop}{0} \setcounter{rem}{0}
\setcounter{tho}{0}
We start by proving the existence result of Theorem \ref{existencer}, for solutions corresponding to a more regular datum. \vskip0.2cm
{\bf {Proof of Theorem \ref{existencer}}} - For any fixed $\mu>0$, denote by $v^\mu$ the solution of \eqref{PFl} corresponding to the initial data
$v_\circ\in W_0^{1,2}(\O)\cap L^\infty(\O)$, satisfying    the
bounds {\sl i), ii)} in Proposition\,\ref{existencemu},  uniformly in  $\mu>0$. Therefore, all the convergence arguments of Proposition\,\ref{existencemu} can be reproduced here in the limit as $\mu\to 0^+$. Let us denote by $u$ the limit function. {\it Mutatis mutandis}, one has only to take care of the convergence of the following term, for any test function $\psi\in C_0^\infty([0,T)\times \O)$: $$
I_4^\mu:=\intll0t(J_\mu(v^\mu)\cdot \nabla v^\mu, \psi)d\tau. $$ 
Actually, we want to show that
$$I_4^\mu \to \intll0t (u\cdot \nabla u, \psi)d\tau, \mbox{ as } \mu\to 0.$$
We write
$$\ba{ll}\dy \vs1 I_4^\mu - \intll0t (u\cdot \nabla u, \psi)d\tau=\intll0t(J_\mu(v^\mu-u)\cdot \nabla v^\mu, \psi)d\tau\\ \dy + \intll0t((J_\mu(u)-u)\cdot \nabla v^\mu, \psi)d\tau
+ \intll0t(u\cdot (\nabla v^\mu-\nabla u), \psi)d\tau. \ea$$
  All the three terms on the right-hand side can easily be shown to go to zero,  by the strong convergence of $v^\mu$ to $u$ and of $\nabla v^\mu$ to $\nabla u$, both in $L^2(\O_T)$, the uniform bound of $\nabla v^\mu$ in $L^2(\O_T)$ and well known convergence properties of mollifiers. \par
  Estimates  {\sl i), ii)} are obtained by lower semicontinuity, while for estimate \eqref{mmtp1r} we argue as in Proposition \ref{existencemu}. Actually, 
  by the strong convergence of the sequence $\{u^{\mu}\}$ in $L^2(0,T; L^2(\O))$ and the boundedness of $ u^{\mu}_t$ uniformly in $\mu$ in $L^2(0,T; L^2(\O))$, one shows that for any fixed $t\geq 0$ the sequence $\{u^\mu(t)\}$ is a Cauchy sequence in $L^2(\O)$ (see \eqref{ado}). Therefore, for any fixed $t\geq 0$, it  strongly converges in $L^2(\O)$, whence a.e. in $\O$. Writing, for any $t\geq 0$ and a.e. in $\O$, 
$$|u(t,x)|\leq |v^\mu(t,x)-u(t,x)|+\|v^\mu(t)\|_\infty\leq |v^\mu(t,x)-u(t,x)|+\|v_\circ\|_\infty,$$
where, in the last step, we have used \eqref{mpvmu} on $v^\mu$,  
passing to the limit on $\mu$, then taking the essential supremum, we also get the validity of \eqref{mmtp1r} for $u$. 
  \vskip0.1cm   
  \begin{coro}\label{existencenomup}{\sl
Under the assumptions of Theorem \ref{existencer} there exists a constant $K$, independent of $\|\nabla v_\circ\|_2$ but depending on $\|v_\circ\|_\infty$,  such that \begin{itemize}
 \item[ii')] $\|t^{\frac{\alpha}{4-p}}\nabla v\|_{L^{\infty}(0,T;L^2(\O))}+\|t^{\frac{\alpha}{2}}D^2 v\|_{L^{2}(0,T; L^\frac{4}{4-p}(\O))}+\|t^\frac 12\,v_t\|_{L^2(\O_T)}\leq K,$
\end{itemize}
\par where $\alpha>\frac{4-p}{p}$. }
\end{coro}
\vskip0.2cm
{\bf {Proof of Theorem \ref{uni}}} - 
Let $u$ and $v$ be two solutions of \eqref{PF}, with the regularity stated in Theorem \ref{existencer}, and let $w:=u-v$ be their difference.  
Differentiating  (\ref{testsmoothst}) with respect to $t$, we find
\be\label{diffU}
 (u_t, \psi)+ (|\nabla u|^{p-2}\nabla u, \nabla\psi)+ (u\cdot\nabla u, \psi)=0,\ee
 for all $\psi \in W^{1,2}(0,T;L^2(\O))\cap L^p(0,T;W_0^{1,p}(\O))\cap L^{p'}((0, T)\times \O),$
 and the same identity holds for $v$ in place of $u$.
We can now take the difference of the above weak formulations written for $u$ and for $v$, and use $w$ as test function. Then, integrating by parts,  we find
\be\label{diffU1}
\begin{array}{ll}\vs1\dy
\frac{1}{2} \frac{d}{dt} ||w||_2^2 + \int_\O (|\nabla u|^{p-2}\nabla u-|\nabla v|^{p-2}\nabla v)\cdot \nabla w dx\\\vs1\hfill \dy 
=-\int_\Omega \left(w \cdot \nabla  u \right)\cdot w \, dx -\int_\Omega \left( v \cdot \nabla w \right)\cdot w \, dx\\ \hfill\dy
=
-\int_\Omega \left(w \cdot \nabla  u \right)\cdot w \, dx +\frac 12\int_\Omega ( \nabla\cdot v) w\cdot  w \, dx
. 
\end{array}
\ee
For the integral on the left-hand side we use the following well known estimate (see, for instance, \cite{DER})
$$(|\nabla u|^{p-2}\nabla u-|\nabla v|^{p-2}\nabla v)\cdot \nabla w\geq C|\nabla w|^2(|\nabla u|+|\nabla v|)^{p-2}.
$$
We apply H\o lder's inequality with exponent $\frac 2q$, $q\leq \frac{4}{4-p}$, and we easily find
$$\ba{ll}\dy\vs1
\|\nabla w\|_q^q= \int_\O \left((|\nabla u|+|\nabla v|)^{p-2}|\nabla w|^2\right)^\frac q2(|\nabla u|+|\nabla v|)^\frac{(2-p)q}{2}\,dx\\\vs1 \dy\leq 
\left(\int_\O (|\nabla u|+|\nabla v|)^{p-2}|\nabla w|^2dx \right)^\frac q2\left(\int_\O(|\nabla u|+|\nabla v|)^\frac{(2-p)q}{2-q}\,dx\right)^\frac{2-q}{2}
\\ \dy \leq c
(\int_\O (|\nabla u|^{p-2}\nabla u-|\nabla v|^{p-2}\nabla v)\cdot \nabla w dx)^\frac q2
\left(\|\nabla u\|_{\frac{(2-p)q}{2-q}}^\frac{(2-p)q}{2}+\|\nabla v\|_{\frac{(2-p)q}{2-q}}^\frac{(2-p)q}{2}\right).\ea$$
Our choice of $q$ ensures that a.e. in $t\in (0,T)$, $$
\left(\|\nabla u(t)\|_{\frac{(2-p)q}{2-q}}^\frac{(2-p)q}{2}+\|\nabla v(t)\|_{\frac{(2-p)q}{2-q}}^\frac{(2-p)q}{2}\right)^\frac 2q\leq 
\|\nabla u\|_{L^\infty(0,T; L^2(\O))}^{2-p}+\|\nabla v\|_{L^\infty(0,T; L^2(\O))}^{2-p}=:\ov C.$$
Using the above estimates in \eqref{diffU1} gives
\be\label{diffU3}
\frac{1}{2} \frac{d}{dt} ||w||_2^2 + \frac{1}{\ov C}||\nabla w||_q^2 \leq \int_\Omega |\left(w \cdot \nabla  u \right)\cdot w \,| dx +\int_\Omega |( \nabla\cdot v) w\cdot  w |\, dx,
\ \mbox{ a.e. in } t\in (0,T).
\ee
We observe that the integrals on the right-hand side are of the same type. Therefore, we just treat one of them. Let us further assume that $q>\frac{12}{7}$. Then, firstly we apply H\o lder's inequality with exponent  $2$, then Gagliardo-Nirenberg's inequality, and, finally, Young's inequality, and  we find
$$\ba{ll}\vs1\dy \int_\Omega |\left(w \cdot \nabla  u \right)\cdot w \,| dx \leq \|\nabla u\|_2 \|w\|_{4}^2\leq 
c\|\nabla u\|_2 \|\nabla w\|_{q}^{2a}\|w\|_{2}^{2(1-a)}\\ \hfill \dy\leq
\eta \|\nabla w\|_{q}^{2}+c(\eta) \|\nabla u\|_2^\frac{1}{1-a}\|w\|_{2}^{2},\ea$$
for any $\eta>0$, with $a:=\frac{3q}{2(5q-6)}.$
Inserting the previous estimates in \eqref{diffU3}, choosing $\eta= \frac{1}{\ov C}$, we get \be\label{diffU4}
\frac{1}{2} \frac{d}{dt} ||w||_2^2 + \frac{1}{2\ov C}||\nabla w||_q^2 \leq c \|\nabla u\|_2^\frac{1}{1-a}\|w\|_{2}^{2},
\ \mbox{ a.e. in } t\in (0,T).
\ee
The application of Gronwall's inequality implies the uniqueness.\chiu
  \vskip0.1cm\begin{lemma}\label{stc}
{\sl Assume that  $\{u^m\}$ is a sequence strongly converging  to a function $u$ in $L^2(0,T; L^2(\O))$ and that the sequence $\{t^\frac 12 u^m_t\}$ is bounded uniformly in $m$ in $L^2(0,T; L^2(\O))$. Then for any $t\geq 0$ the sequence $\{u^m(t)\}$ strongly converges to $u(t)$ in $L^2(\O)$. }
\end{lemma}
\Pr By elementary calculations and H\o lder's inequality we get
$$\ba{ll}\vs1\dy t\|u^{m_1}(t)-u^{m_2}(t)\|_2^2=\int_0^t\frac{d}{d\tau}(\tau\|u^{m_1}(\tau)-u^{m_2}(\tau)\|_2^2)d\tau\\\vs1\dy =\int_0^t\|u^{m_1}(\tau)-u^{m_2}(\tau)\|_2d\tau+
2
\int_0^t\tau\,(u^{m_1}(\tau)-u^{m_2}(\tau))\cdot (u^{m_1}_\tau(\tau)-u_\tau^{m_2}(\tau))d\tau\\ \vs1\dy
\vs1 \leq\int_0^t\|u^{m_1}(\tau)-u^{m_2}(\tau)\|_2d\tau+
2
\int_0^t\tau^\frac 12\,\|u^{m_1}(\tau)-u^{m_2}(\tau)\|_2 \tau^\frac 12\,\|u^{m_1}_\tau(\tau)-u_\tau^{m_2}(\tau)\|_2d\tau
\\\dy\leq
\int_0^t\|u^{m_1}(\tau)-u^{m_2}(\tau)\|_2d\tau\\ \dy \hfill+
2 t^\frac 12(\int_0^t\|u^{m_1}(\tau)-u^{m_2}(\tau)\|_2^2d\tau)^\frac 12(\int_0^t(\tau^\frac 12\|u_\tau^{m_1}(\tau)-u_\tau^{m_2}(\tau)\|_2)^2d\tau)^\frac 12.
\ea $$
By the strong convergence of the sequence $\{u^m\}$ in $L^2(0,T; L^2(\O))$ and the assumption of boundedness of $t^\frac 12 u^m_t$ uniformly in $m$ in $L^2(0,T; L^2(\O)$, we find that for any fixed $t\geq 0$ the sequence $\{u^m(t)\}$ is a Cauchy sequence in $L^2(\O)$, hence strongly converges in $L^2(\O)$. \chiu

  \begin{prop}\label{existencenon}{\sl Assume that $u_\circ$
belongs to $L^\infty(\O)$. Then there exists a solution $u$
of system \eqref{PF} in the sense of Definition\,\ref{defnomu}, for any $T>0$.
Moreover $u$ satisfies the following properties:\begin{itemize}
\item[i)]
$\dy \|u\|_{L^{\infty}(0,T;L^2(\O))}+\|\nabla u\|_{L^{p}(\O_T)}\leq K_1(\|u_\circ\|_\infty, T)\,;$
 \item[ii)]  
  $t^{\frac{\alpha}{4-p}}\nabla u\in L^{\infty}(0,T;L^2(\O))$, $t^{\frac{\alpha}{2}}D^2 u\in L^{2}(0,T; L^\frac{4}{4-p}(\O)),$  $t^\frac 12\,u_t\in L^2(\O_T);$
\item[iii)]$\dy\lim_{t\to 0^+}\|u(t)-u_\circ\|_2=0.$
\end{itemize}
\par where $\alpha>\frac{4-p}{p}$. }
\end{prop}
  \Pr Let $\dy\{ u_\circ^m(x)\}$ be a sequence
in $C_0^\infty(\O)$ strongly converging to $u_{\circ}$ in $L^2(\O)$, such that $\|u_\circ^m\|_\infty\leq \|u_\circ\|_\infty$, 
and let $\{u^{m}\}$ be the sequence of solutions of system
(\ref{PF}) corresponding to the initial data $\{u_\circ^m(x)\}$.
From Corollary\,\ref{existencenomup}, 
 for any data $u_\circ^m$, there exists a solution $u^m$
 of \eqref{PF} that satisfies {\sl i)} in Theorem  \ref{existencer} and such that $t^{\frac{\alpha}{2}}\,u^m
 \in L^2(0,T; W^{2,\frac{4}{4-p}}(\O))$ and $t\,u^m_{t} \in L^{2}(\O_T)$.
 Then,  as $m\to \infty$, the convergence results \eqref{wnoee}-\eqref{wnoe}, in Proposition \ref{existencemu}, hold, while the convergence results \eqref{weg}--\eqref{wedt} and their consequences, \eqref{dn1} and \eqref{dn2}, hold in a time interval $(\ve, T)$, for 
any $\ve>0$, by Corollary \ref{existencenomup}. \par   
Further, we observe that the sequence $\{u^m\}$ is a Cauchy sequence in $L^p(0,T; L^2(\O))$, due to the weak convergence of $\{u^m(t)\}$ to $u(t)$ in $L^2(\O)$, uniformly in $t\geq 0$, and the uniform boundedness of $\{\nabla u^m(t)\}$ in $L^p(\O_T)$. Indeed, the generalized Friedrich's lemma (see for instance \cite{LSU}, Lemma II.2.4) ensures that, denoting by 
 $(\psi_k)$ an orthogonal basis of $L^2(\O)$, for any $\eta>0$, there exists a number $N_\eta$ such that
$$\intl0^T\!\|u^{m_1}(t)-u^{m_2}(t)\|_2^p dt\leq 
\!\intl0^T\!\!\big(\sum_{k=1}^{N_\eta}(u^{m_1}(t)-u^{m_2}(t), \psi_k)^2\big)^\frac p2 dt
+\eta\!\intl0^T\!\|\nabla u^{m_1}(t)-\nabla u^{m_2}(t)\|_p^p dt.$$
Therefore the sequence $\{u^m\}$ strongly converges to $u$ in $L^p(0,T; L^2(\O))$. This convergence and the uniform boundedness of $\{u^m\}$ in $L^\infty(0,T; L^2(\O))$, ensure the strong convergence in $L^2(\O_T)$, by interpolation. 
Then, the sequence $\{u^m\}$ satisfies the hypotheses of Lemma \ref{stc}, and therefore for any $t\geq 0$, $\{u^m(t)\}$ strongly converges to $u(t)$ in $L^2(\O)$, whence a.e. in $\O$. 
For any $t\geq 0$ and a.e. in $\O$, 
$$|u(t,x)|\leq |u^m(t,x)-u(t,x)|+\|u^m(t)\|_\infty\leq 
|u^m(t,x)-u(t,x)|+\|u_\circ\|_\infty,$$
where in the last step we have used the validity of estimate \eqref{mmtp1r} on the sequence $\{u^m(t)\}$ and the assumption $\|u_\circ^m\|_\infty\leq \|u_\circ\|_\infty$. Passing to the limit on $m$ and then taking the essential supremum, we get \eqref{mmtp1}. 
\par We observe that, since $u^m(t)\in C([0,T);L^2(\O))$, we can use Remark \ref{testsmoothN} 
 and obtain that, for any $\ve>0$, 
$$\ba{ll}\dy\vs1 \int_\ve^t
\left[(u^m,\psi_\tau)-\left(|\nabla u^m|^{p-2}\,\nabla
u^m,\nabla \psi\right)-(u^m\cdot \nabla u^m, \psi)\right]d\tau\\ \hskip3.5cm \dy =(u^m(t),\psi(t))+(u^m(\ve),
\psi(\ve)), \ \mbox{ for any } \ \psi\in C_0^\infty([0,T)\times \O). 
\ea$$
 Passing to the limit as $m\to \infty$, we can argue as before and obtain
 $$\ba{ll}\dy\vs1 \int_\ve^t
\left[(u,\psi_\tau)-\left(|\nabla u|^{p-2}\,\nabla
u,\nabla \psi\right)-(u\cdot \nabla u, \psi)\right]d\tau\\ \hskip4cm \dy =(u(t),\psi(t))+(u(\ve),
\psi(\ve)),\ 
 \mbox{for all }\psi \in C_0^\infty([0,T)\times \O).\ea$$ 
We remark that, by standard arguments, one can verify that $u$ is weakly continuous in $L^2(\O)$. Since $$|(u(\ve),
\psi(\ve))-(u(0),
\psi(0))|\leq |(u(\ve),
\psi(\ve)-\psi(0))|+|(u(\ve)-u(0),
\psi(0))|, $$
 from the $L^2$-weak continuity of $u$ and the regularity of $\psi$ one has
 $$\lim_{\ve\to 0^+}|(u(\ve),
\psi(\ve))-(u(0),
\psi(0))|=0.$$ From this remark and the properties of the limit function, i.e. $u\in L^\infty(0,T; L^2(\O))\cap L^p(0,T; W_0^{1,p}(\O))$, using
 the continuity of the Lebesgue integral, we can pass to the limit as $\ve\to 0^+$, and find that $u$ satisfies the integral formulation of Definition \ref{defnomu}.\par As $u\in L^{p}(0,T; W_0^{1,p}(\O))$ and $u_t\in L^{2}(\ve,T; L^2(\O))$,  the $L^2(\O)$-continuity property of the solution $u$ in $[\ve, T)$  follows. Therefore, using the energy estimate and the $L^2$-weak continuity, one gets $$\dy\lim_{t\to
0^+}\|u(t)-u_\circ\|_2^2=\dy\lim_{t\to
0^+}(\|u(t)\|_2^2+\|u_\circ\|_2^2-2(u(t), u_\circ))=0
\,,$$ which, together with $u\in C([\ve,T); L^2(\O))$ for any $\ve>0$, gives $u\in C([0,T); L^2(\O))$. The continuity in $L^q(\O)$, for any $q\in [1,\infty)$ follows by interpolation. \chiu 

\begin{rem}\label{testsmooth}
As already pointed out in Remark \ref{testsmooth}, we observe that, due to regularities obtained for $u$, for any $\ve>0$ if we write the weak formulation \eqref{testsmoothst} in $(\ve, t)$, we can take the test functions 
$\psi$ in the space $W^{1,2}(\ve,T;L^2(\O))\cap L^p(\ve,T;W_0^{1,p}(\O))\cap L^{p'}(\ve,T; L^{p'}(\O_T))$. 
\end{rem}

 \section*{\normalsize 5. Proof of the extinction theorem}
\renewcommand{\theequation}{5.\arabic{equation}}
\renewcommand{\thetho}{5.\arabic{tho}}
\renewcommand{\thedefi}{5.\arabic{defi}}
\renewcommand{\therem}{5.\arabic{rem}}
\renewcommand{\theprop}{5.\arabic{prop}}
\renewcommand{\thelemma}{5.\arabic{lemma}}
\renewcommand{\thecoro}{5.\arabic{coro}}
\setcounter{equation}{0} \setcounter{coro}{0} \setcounter{lemma}{0}
\setcounter{defi}{0}
%
 We prove the result by considering system \eqref{PF}  with any $\delta>0$. Let us consider the solution $u$ given in Theorem \ref{existence}, for which both the integral formulation in Definition \ref{defnomu} and in Remark \ref{testsmooth} hold.
Differentiating  (\ref{testsmoothst}) with respect to $t$, we find
\be\label{diff1}
 (u_t, \psi)+ (|\nabla u|^{p-2}\nabla u, \nabla\psi)+\delta (u\cdot\nabla u, \psi)=0,\ee
 for all $\psi \in W^{1,2}(\ve,T;L^2(\O))\cap L^p(\ve,T;W_0^{1,p}(\O))\cap L^{p'}((\ve, T)\times \O).$
We can now use $u$ as test function in \eqref{diff1}, then apply Hölder and Young inequalities, and find
$$
\frac{1}{2} \frac{d}{dt} ||u||_2^2 + ||\nabla u||_p^p \leq \delta \left|\int_\Omega \left( u \cdot \nabla  u \right)\cdot u \, dx \right| \leq \frac 1p ||\nabla u||_p^p + \frac{\delta}{p'} ||u^2||_{p'}^{p'},\ \forall t\in (\ve,T).
$$
So, by convexity inequality,
$$
\frac{1}{2} \frac{d}{dt} ||u||_2^2 + \frac 1{p'} ||\nabla u||_p^p \leq \frac{\delta}{p'} ||u||_{2p'}^{2p'} \leq \frac{\delta}{p'}\left( ||u||_\infty^a ||u||_2^{1-a}\right)^{2p'},
$$
where $\frac{1}{2p'} = \frac{1-a}{2}$.
Since, by the Sobolev embedding Theorem, it is known that
$$
\begin{array}{ccc}
||\nabla u||_p^p \geq c_s||u||_{p^\star}^p \geq \gamma ||u||_2^p
\end{array}$$
with $\gamma=\gamma(\Omega, p)$,  we get
$$
\begin{array}{ccc}\dy
\frac{1}{2} \frac{d}{dt} ||u||_2^2 + \frac{\gamma}{p'} || u||_2^p \leq \frac{\delta}{p'}||u||_\infty^{\frac{2}{p-1}} ||u||_2^{2},\ \forall t\in (\ve,T),
\end{array}$$
that implies
\begin{equation}
\begin{array}{ccc}
\label{ineq}\dy
\frac{d}{dt} ||u||_2 \leq \frac{1}{p'} ||u||_2^{p-1} \left( {\delta} ||u||_2^{2-p}||u_\circ||_\infty^{\frac{2}{p-1}} - \gamma\right), \forall t\in (\ve,T).
\end{array}
\end{equation}
Setting $\varphi(t):= ||u(t)||_2$, we want to prove that $\varphi ' (t) \leq 0$ for all \mbox{$t \in (\ve,T)$}.
Due to the assumptions on the initial datum, there exists a $\ov \sigma$ such that for any $\sigma\in (0, \ov \sigma)$
$$\delta||u_\circ||_2^{2-p}||u_\circ||^{\frac{2}{p-1}}_\infty(1+\sigma)<{\gamma}.$$ From the strong continuity of $u$ in $L^2(\O)$, there exists $\ov \ve$ such that 
there holds
\be\label{starC}||u(\ve)||_2^{2-p}\leq (1+\sigma)||u_\circ||_2^{2-p}, \ \forall \ve\in (0,\ov\ve].\ee
This, together with \eqref{ineq}, implies that  $\varphi'(\ov \ve)<0$ and $\varphi(t)<\varphi(\ov \ve)$ in a right neighborhood of $\ov\ve$.  By \eqref{ineq}, this ensures that $\varphi$ is decreasing in the same neighborhood of $\ov\ve$.  On the other hand one readily proves that $\varphi$ decreases in the whole interval $(\ov\ve, T)$,  reasoning by contradiction. Indeed, let $\ov t:=\sup\{s:\varphi(\tau)<\vp(\ov \ve), \forall \tau\in (\ov \ve, s)\}$. By continuity $\varphi(\ov t)=\varphi(\ov \ve)$ and, by \eqref{ineq},   $\varphi'(\ov t)<0$, that gives the contradiction.   \par
Now, in order to conclude, we restart from (\ref{ineq}) and, using the monotonicity of $\varphi(t)$, we get
\be\label{nss}
\begin{array}{ccc}
\varphi '(t) \leq \frac{1}{p'}\varphi ^{p-1}(t) \left({\delta} ||u_0||_\infty^{\frac{2}{p-1}} \|u_\circ\|_2^{2-p}(1+\sigma)- \gamma \right) , \forall t\in (\ve, T).
\end{array}
\ee
Since $\varphi'(t)\leq 0$, for any $t\geq \ov \varepsilon$, if there exists $t_0$ such that $\varphi(t_0)=0$, then $\varphi(t)\equiv 0$, for any $t\geq t_0$. Hence, we suppose $\varphi(\ov\ve)>0$ and we divide \eqref{nss} by $\varphi^{p-1}(t)$ as long as $\varphi(t)>0$. Integrating in time, there holds
$$\begin{array}{ccc}
\varphi^{2-p} (t) \leq \varphi^{2-p} (\ve) + \frac{2-p}{p'}\, t \left( {\delta} ||u_\circ||_\infty^{\frac{2}{p-1}} \|u_\circ\|_2^{2-p}(1+\sigma)- \gamma \right).
\end{array}$$
Since if $$t=\frac{p'(1+\sigma)||u_\circ||_2^{2-p}}{(2-p) \big(\gamma-{\delta} ||u_\circ||_\infty^{\frac{2}{p-1}} (1+\sigma)||u_\circ||_2^{2-p}  \big)   }:=T^*, $$ by \eqref{starC} the right-hand side is non-positive, it follows that there exists $t_0\leq T^*$ such that $\vp(t_0)=0$, hence
$$\varphi(t) = ||u(t)||_2 \equiv 0 \text{ for all } t \geq t_0.$$
Passing to the limit as $\sigma\to 0^+$ in $T^*$, we obtain the result.
\chiu
\begin{rem}
As $\delta \to 0^+$, $T^* \to \frac{p'||v_0||^{2-p}_2}{(2-p)\gamma } $ according to \cite{DB}. Actually, if we decide not to carry out a precise comparison with the constant in \cite{DB}, the proof of the previous theorem can be considerably simplified. For example, we can avoid introducing the additional constant $\sigma$, by simply imposing $$2\delta||u_\circ||_2^{2-p}||u_\circ||^{\frac{2}{p-1}}_\infty<{\gamma}.$$ .\par Further, as $p \to 2^-$, $T^* \to + \infty$. This is in agreement with what is known for the Navier-Stokes equations, for which the phenomenon of extinction in a finite time does not exist, see \cite{GaS}. 
\end{rem}
\vskip0.1cm
 {\bf Acknowledgment} -
The paper is
 performed under the
auspices of GNFM-INdAM.

\end{document}